\documentclass[article,hidelinks,onefignum,onetabnum]{siamart250211}

\usepackage{algorithm}
\usepackage{algpseudocode}
\usepackage{amsfonts,amsmath,amssymb}
\usepackage{booktabs}
\usepackage{epstopdf}
\usepackage{graphicx}
\usepackage{mathtools}
\usepackage{pgf, tikz}
\usepackage{pgfplots}
\usepackage{pgfplotstable}
\usepackage{stmaryrd}
\usepackage{xcolor, soul}
\pgfplotsset{compat=1.14}
\usetikzlibrary{shapes.geometric, arrows,backgrounds,fit,positioning}
\ifpdf
  \DeclareGraphicsExtensions{.eps,.pdf,.png,.jpg}
\else
  \DeclareGraphicsExtensions{.eps}
\fi

\pgfplotsset{
tick label style={font=\tiny},
legend style={font=\scriptsize, draw=none},
xlabel style={yshift=+0.5ex},
ylabel style={yshift=-1.0ex}
}

\newcommand{\R}{\mathbb R}
\newcommand{\N}{\mathbb N}
\newcommand{\e}{\mathrm{e}}
\newcommand{\T}{\mathcal{T}}
\newcommand{\tol}{\mlvar{\mathrm{tol}}}
\newcommand\abrev[1]{{\color{black}#1}}
\newcommand\rev[1]{{\color{black}#1}}
\DeclarePairedDelimiter{\norm}{\lVert}{\rVert}  
\newcommand{\scal}[2]{\left(#1,#2\right)}

\mathchardef\standardl=\mathcode`l
\mathcode`l=\ell
\newcommand{\deactivatel}{\mathcode`l=\standardl}
\makeatletter
\edef\operator@font{\operator@font\noexpand\deactivatel}
\makeatother

\newcommand{\mlvar}[1]{\mathit{\deactivatel #1}}

\newsiamremark{remark}{Remark}
\newsiamremark{hypothesis}{Hypothesis}
\crefname{hypothesis}{Hypothesis}{Hypotheses}
\newsiamthm{claim}{Claim}
\newsiamremark{fact}{Fact}
\crefname{fact}{Fact}{Facts}

\tikzset{block/.style={rectangle, draw, text width=4.2em, text centered, minimum height=3em, rounded corners=5pt},
title/.style={font=\bfseries},
line/.style={draw, -latex'},
int/.style={draw, minimum size=2em, font=\upshape}}

\headers{Adaptive multimesh rational approximation scheme}{A. Bespalov, R. Bulle}

\ifpdf
\hypersetup{
  pdftitle={An adaptive multimesh rational approximation scheme for the spectral fractional Laplacian},
  pdfauthor={A. Bespalov, R. Bulle}
}
\fi

\title{An adaptive multimesh rational approximation scheme for the spectral fractional Laplacian}

\author{Alex Bespalov\thanks{School of Mathematics, University of Birmingham, Edgbaston, Birmingham B15 2TT, UK (\email{a.bespalov@bham.ac.uk})}
\and Rapha{\"e}l Bulle\thanks{Team-project MIMESIS, Inria de l'Universit{\' e} de Lorraine, MLMS team, Icube, Universit{\' e} de Strasbourg, CNRS UMR 7357,
2 Rue Marie Hamm, 67000 Strasbourg, France (\email{raphael.bulle@inria.fr})}}

\begin{document}

\maketitle

\begin{abstract}
We propose a novel \emph{multimesh} rational approximation scheme for the numerical
solution of the (homogeneous) Dirichlet problem for the spectral fractional Laplacian.
The scheme combines a rational approximation of the function $\lambda \mapsto \lambda^{-s}$
with a family of finite element discretizations of parameter-dependent, non-fractional
partial differential equations (PDEs).
The key idea that underpins the proposed scheme is that each parametric PDE
is numerically solved on an individually tailored finite element mesh.
This is in contrast to existing \emph{single-mesh} approaches that employ the same
finite element mesh across all parametric PDEs.
We develop an a posteriori error estimation strategy for the proposed
rational approximation scheme and design
an adaptive multimesh refinement algorithm.
Numerical experiments demonstrate that our adaptive multimesh approach achieves faster convergence rates
than uniform mesh refinement and yields
significant reductions in computational costs---both in terms of the overall number of degrees of freedom
and the actual runtime---when compared to the corresponding adaptive algorithm in a single-mesh setting.
\end{abstract}



\section{Introduction}

Fractional partial differential equations (FPDEs) arise in a wide range of applications
including anomalous diffusion, 
porous media, 
phase separation in fluids, 
discontinuous deformations (fractures),
and spatial statistics 
(see~\cite{2018_bonito, DEliaDGGTZ_20_NMN, 2020_lischke} and the references therein).
They have become a powerful mathematical tool in these application areas due to their ability to model nonlocal
phenomena and processes that are characterized by interactions at a distance.
The development of effective bespoke numerical methods and algorithms for FPDEs is also of significant interest,
as the standard computational techniques applied to them generally incur much greater computational
costs than if applied to classical (non-fractional) PDE-based models.

In this work, we consider one particular fractional derivative operator,
the fractional Laplacian $(-\Delta)^s$ ($0 < s < 1$), that has received the most attention from both
the mathematical analysis and numerical analysis communities.
Among many non-equivalent definitions of the fractional Laplacian on a bounded domain $\Omega$,
our focus here is on the \emph{spectral} fractional Laplacian
(we give \rev{a} precise definition in section~\ref{sec:spectral:FL}).
For a given $s \in (0,1)$ and 
$f \in L^2(\Omega)$, we consider the following
boundary value problem for the spectral fractional Laplacian:
find $u : \Omega \to \R$ such that
\begin{equation}
    (-\Delta)^s u = f\quad \text{in } \Omega,\qquad u = 0\quad \text{on } \partial \Omega.
    \label{eq:fractional_laplacian_strong}
\end{equation}
%
This problem can be discretized directly using, e.g., the finite element method (FEM).
However, due to the nonlocal nature of the fractional Laplacian, such a discretization would lead to a dense linear system,
see, e.g.,~\cite{2020_lischke}.
The resulting linear system can become computationally intractable if the number of degrees of freedom is very large, which is
often the case when solving 3D problems and/or using a very fine computational mesh.
In addition, for small values of $s$, the solution to \eqref{eq:fractional_laplacian_strong} exhibits boundary layers,
calling for approximations that employ adaptive and/or anisotropic mesh refinement,
see, e.g.,~\cite{2015_chen, 2023_banjai}.

Several alternative discretization strategies for 
problems involving the fractional Laplacian (and applicable to more general nonlocal problems) are available in the literature.
Among these strategies we mention
  the methods that employ the Dirichlet-to-Neumann maps to reformulate~\eqref{eq:fractional_laplacian_strong} in 2D
  as a three-dimensional non-fractional problem (see~\cite{2015_nochettob, BanjaiMNOSS_19_TFS}),
  reduced basis methods that construct a discrete space of small dimension where
  the solution $u$ is sought (see~\cite{2019_antil}),
  the semi-group method that transforms the elliptic problem~\eqref{eq:fractional_laplacian_strong}
  into a non-fractional parabolic problem~(see~\cite{2018_cusimano}),
  and
  the methods that employ rational approximations of the function $\lambda \mapsto \lambda^{-s}$
  in order to reformulate a fractional PDE as a set of independent non-fractional
  parameter-dependent PDEs
  (see~\cite{2015_bonito,2016_bonito,2020_harizanov,2020_vabishchevich,2023_bulle}).
%
In particular, the strategies combining rational approximations of $\lambda^{-s}$ and
finite element discretizations of non-fractional parametric PDEs have become popular.
In these strategies, called \emph{rational approximation schemes}, instead of assembling and solving one \emph{dense} linear system,
the FEM is used to assemble several independent \emph{sparse} linear systems that can be solved in parallel.
The approximate solution to problem~\eqref{eq:fractional_laplacian_strong}
is then obtained as a linear combination of those independently generated finite element approximations.
The existing rational approximation schemes employ the same finite element space
(associated with a single finite element mesh on $\Omega$) across \emph{all} non-fractional parametric PDE problems;
we refer to this approach as a \emph{single-mesh} scheme.

A posteriori error estimation and adaptive solution techniques are fairly recent topics in the numerical analysis of FPDEs.
The following a posteriori error estimation strategies have been developed
to guide adaptive mesh refinement in the FEM-based methods for fractional problems:
residual-based error estimators~\cite{2017_ainsworth, 2021_faustmann},
a gradient recovery based error estimator~\cite{2017_zhao},
and hierarchical error estimators~\cite{2015_chen, 2023_bulle}.

In this work, 
we propose a novel \emph{multimesh framework} for rational approximation schemes where
the underlying meshes (and hence FEM approximations) can
be different for different non-fractional parametric problems.
Within this framework, we develop an a posteriori error estimation strategy that drives \rev{the} adaptive multimesh refinement algorithm for the numerical solution of problem~\eqref{eq:fractional_laplacian_strong}.
We emphasize that our multimesh construction, the error estimation strategy, and the adaptive algorithm
are generic, in the sense that they can be used with any rational approximation of the function $\lambda^{-s}$
(see~\S\ref{sec:rational:approx} below).
Our practical implementation, however, 
employs a specific rational approximation---the one proposed by Bonito and Pasciak in~\cite{2015_bonito}---but
can be easily adapted to other approximations such as, e.g., the BURA-based schemes~\cite{2020_harizanov}.
Numerical experiments demonstrate two main advantages of our adaptive multimesh framework:
faster convergence rates (particularly, for smaller values of $s \in (0,1)$) 
than the rates for uniform mesh refinement (see~\cite{2015_bonito}) and
significant reductions in computational costs---measured by the overall number of degrees of freedom
and the actual computational time---when compared to the
adaptive algorithm proposed in~\cite{2023_bulle} in a single-mesh setting.

The paper is organized as follows.
In section~\ref{sec:spectral:FL}, we define the spectral fractional Laplacian operator
$(-\Delta)^s$ in~\eqref{eq:fractional_laplacian_strong}.
Section~\ref{sec:multimesh:scheme} introduces a multimesh rational approximation scheme
for problem~\eqref{eq:fractional_laplacian_strong},
including the underlying rational approximation of the function $\lambda^{-s}$,
finite element discretizations as well as the concept of
union mesh that plays a key role in computing the fully discrete solution and
in a posteriori error analysis.
In section~\ref{sec:error:estimation}, we discuss local hierarchical error indicators and
propose two a posteriori estimates of the overall discretization error in the $L^2$-norm.
An adaptive multimesh refinement algorithm is formulated in section~\ref{sec:adaptive:refinement}.
The effectiveness of two error estimation strategies and the performance of the proposed
adaptive algorithm 
are assessed in a series of numerical experiments presented in
section~\ref{sec:numerical_results}.

\section{The spectral fractional Laplacian} \label{sec:spectral:FL}

\abrev{Let $\Omega \subset \R^d$ ($d=2, 3$) be a connected bounded domain
with polygonal/poly\-hedral boundary $\partial \Omega$ and
let $\omega$} be any open connected subset of $\overline \Omega$.
We denote by $L^2(\omega)$ the space of square integrable functions over $\omega$ \abrev{with
the usual inner product $\scal{\cdot}{\cdot}_{2,\omega}$ and the associated norm $\norm{\cdot}_{2,\omega}$}.
We denote by $H^1(\omega)$ the Sobolev space of functions with \abrev{first-order} weak derivatives in $L^2(\omega)$;
\abrev{it is} endowed with \abrev{the} inner product
$\scal{\nabla\cdot}{\nabla \cdot}_{2,\omega} + \scal{\cdot}{\cdot}_{2,\omega}$\abrev{, and the associated squared norm is given by}
$\norm{\nabla \cdot}_{2,\omega}^2 + \norm{\cdot}_{2,\omega}^2$.
\abrev{The subspace of functions in $H^1(\Omega)$ with zero trace on the boundary $\partial\Omega$ is denoted by $H^1_0(\Omega)$.
The inner product and the norm in $H^1_0(\Omega)$ are given by
$\scal{\nabla \cdot}{\nabla \cdot}_{2,\Omega}$ and $\norm{\nabla \cdot}_{2,\Omega}$, respectively.}
\abrev{In what follows, when $\omega = \Omega$, we will}
omit the dependence \abrev{on} $\omega$ in the subscripts \abrev{of inner products and norms}.

Let $\abrev{\{(\psi_i, \lambda_i);\ i \in \N\}} \subset L^2(\Omega) \times \R_{>0}$
be the spectrum of the \abrev{standard negative Laplacian} $-\Delta$ on $\Omega$
\abrev{with zero} Dirichlet boundary condition on $\partial \Omega$.
\abrev{In other words, the eigenvalues $\lambda_i \in \R_{>0}$ and the eigenfunctions $\psi_i: \Omega \to \R$}
are defined by the following eigenvalue \abrev{problem:}
\begin{equation*}
    -\Delta \psi_i = \lambda_i\psi_i\quad \text{in }\Omega, \quad \psi_i = 0\quad \text{on }\partial \Omega.
\end{equation*}
We assume that the eigenvalues $\lambda_i$ are sorted in increasing order and we denote by $\lambda_0 \in \R_{>0}$ a lower bound such that
\begin{equation*}
  \lambda_0 \leqslant \lambda_1 \leqslant \ldots \leqslant \lambda_i \leqslant
  \lambda_{i+1} \leqslant \ldots.
\end{equation*}
Furthermore, the \abrev{set} $\{\psi_i;\ i\in \N\}$ is an orthonormal basis in $L^2(\Omega)$.

For any $s \in (0, 1)$, \abrev{let us now consider the spectral fractional Laplacian operator $(-\Delta)^{s}$}
with \abrev{zero} Dirichlet boundary condition.
\abrev{This is a pseudo-differential operator that} is defined via its action on the eigenfunctions \abrev{$\psi_i$ of the standard Laplacian as follows}:
\begin{equation*}
  (-\Delta)^{s} \psi_i = \lambda_i^{s} \psi_i\quad \forall i \in \N.
\end{equation*}
For \abrev{any} function
\abrev{$v \in {\mathbb H}^{s}(\Omega) := \Big\{ v \in L^2(\Omega),\ \sum\limits_{i=1}^{\infty} \lambda_i^{s} \scal{v}{\psi_i}_2^2 < \infty \Big\}$},
we have
\begin{equation*}
  (-\Delta)^{s} v := \sum_{i=1}^{\infty} \lambda_i^{s} \scal{v}{\psi_i}_2 \psi_i.
\end{equation*}
Using this representation 
we conclude that
the solution $u$ to \abrev{problem}~\eqref{eq:fractional_laplacian_strong} admits the following expansion:
\begin{equation}
    u = \sum_{i=1}^{\infty} \lambda_i^{-s} \scal{f}{\psi_i}_2 \psi_i.
    \label{eq:fractional_solution_spectral_expansion}
\end{equation}

\section{A multimesh rational approximation scheme} \label{sec:multimesh:scheme}

In this section, we describe two components of rational approximations schemes:
a rational approximation of the function $\lambda \mapsto \lambda^{-s}$ and
the finite element discretization.
We also introduce the union mesh as a key ingredient of our \emph{multimesh} construction.

\subsection{Rational approximations}
\label{sec:rational:approx}

\abrev{The idea of rational approximation schemes in the context of solving problem~\eqref{eq:fractional_laplacian_strong}
hinges on approximating} the function $\lambda \mapsto \lambda^{-s}$ for $\lambda \in [\lambda_0,+\infty)$ \abrev{by rational functions}
$\lambda \mapsto \mathcal Q_s^{\kappa}(\lambda)$ of the form
\begin{equation}
  \mathcal Q_s^{\kappa}(\lambda) :=
  C(s,\kappa) \sum_{l=1}^{N(\kappa)} a_l(s,\kappa) \big(c_l(s,\kappa) + b_l(s,\kappa) \lambda\big)^{-1}.
  \label{eq:definition_rational_function}
\end{equation}
\abrev{Here,} $\kappa \in \R_{>0}$ encodes the fineness of \abrev{this approximation},
$\N_{>0} \ni N(\kappa) \to +\infty$ as $\kappa \to 0$,
\abrev{ and} the coefficients
$C,\, a_l,\, b_l,\, c_l \in \R_{>0}$ for all $l=1,\ldots,N(\kappa)$
are chosen so that $\mathcal Q_s^{\kappa}$ converges exponentially fast \abrev{to $\lambda^{-s}$
(uniformly on the interval $[\lambda_0,+\infty)$)} as $\kappa \to 0$.
In other words, there exists $\varepsilon_s(\kappa) \in \R_{>0}$ independent of $\lambda$ such that
\begin{equation}
  \left\lvert \lambda^{-s} - \mathcal Q_s^{\kappa}(\lambda) \right\rvert \leqslant
  \varepsilon_s(\kappa)\quad \forall \lambda \in [\lambda_0,+\infty),
  \label{eq:rational_scheme_convergence}
\end{equation}
where \abrev{$\varepsilon_s(\kappa) \to 0$ exponentially fast as $\kappa \to 0$}.

\abrev{The rational approximation $\mathcal Q_s^{\kappa}$ defined} in~\eqref{eq:definition_rational_function} can be used to derive a semi-discrete
approximation \abrev{of the solution $u$ to problem~\eqref{eq:fractional_laplacian_strong}.
Specifically, replacing} $\lambda_i^{-s}$ with $\mathcal Q_s^{\kappa}(\lambda_i)$ for all $i \in \N$ in
\abrev{representation~\eqref{eq:fractional_solution_spectral_expansion} of the solution $u$, the
following semi-discrete approximation of $u$ is obtained
(see, e.g.,~\cite{2015_bonito} for full details):} 
\begin{equation}
  u_{\kappa} := C(s,\kappa) \sum_{l=1}^{N(\kappa)} a_l(s,\kappa) w_l,
  \label{eq:rational_approximation_definition}
\end{equation}
where each function \abrev{$w_l : \Omega \to \R$}
(\abrev{$l = 1,\ldots,N(\kappa)$}) solves \abrev{the corresponding} \emph{non-fractional}
reaction-diffusion problem
\begin{equation}
  \abrev{-b_l(s,\kappa) \Delta w_l + c_l(s,\kappa) w_l = f\quad
  \text{in } \Omega,\qquad w_l = 0\quad \text{on } \partial \Omega.}
  \label{eq:continuous_parametric_problems}
\end{equation}
\abrev{The fully discrete approximation of $u$ is then obtained by 
discretizing each parametric problem~\eqref{eq:continuous_parametric_problems}
(using, e.g., the FEM, see~\S\ref{subsec:fe_discretization})
and replacing the functions $w_l$ in~\eqref{eq:rational_approximation_definition}
with their (finite element) approximations.
}

\subsubsection{\abrev{The} Bonito--Pasciak rational \abrev{approximation}}

For a fixed $s \in (0,1)$, the Bonito--Pasciak (BP) rational \abrev{approximation} scheme \abrev{for solving problem~\eqref{eq:fractional_laplacian_strong}}
is based on the following \abrev{identity}
\abrev{that is} derived from the Balakrishnan formula \abrev{in}~\cite{1960_balakrishnan}:
\begin{equation}
    \lambda^{-s} = \frac{2\sin(\pi s)}{\pi}
    \int_{-\infty}^{+\infty} \e^{2s y} \left(1 + \e^{2y}
    \lambda\right)^{-1}\, \mathrm{d}y.
    \label{eq:balakrishnan}
\end{equation}
\abrev{Hence, for} a given fineness parameter $\kappa \in \R_{>0}$, the rational \abrev{approximation}
$\mathcal Q_s^{\kappa}(\lambda)$ of $\lambda^{-s}$ is obtained \abrev{in~\cite{2015_bonito}}
by \abrev{discretizing} the integral in~\eqref{eq:balakrishnan}
via a rectangle quadrature rule:
\begin{equation}
    \mathcal Q_s^{\kappa}(\lambda) := \frac{2 \kappa \sin(\pi s)}{\pi
    } \sum_{j=-M_-(\kappa)}^{M_+(\kappa)} \e^{2sj\kappa} \left(1+ \e^{2j\kappa} \lambda\right)^{-1},
    \label{eq:bp_rational_approximation}
\end{equation}
where
\begin{equation*}
    M_-(\kappa) := \left\lceil \frac{\pi^2}{4s\kappa^2} \right\rceil
    \quad\text{and}\quad
    M_+(\kappa) := \left\lceil \frac{\pi^2}{4(1-s)\kappa^2} \right\rceil
\end{equation*}
with $\lceil \cdot \rceil$ \abrev{denoting} the ceiling function.

Setting $l := j + M_-(\kappa) + 1$ and $N(\kappa) := M_+(\kappa) + M_-(\kappa) + 1$, the right-hand side of~\eqref{eq:bp_rational_approximation} can be written in the generic form given by~\eqref{eq:definition_rational_function} with
\begin{subequations}
  \label{eq:bp_rational_approximation_coeffs}
  \begin{align}
  \label{eq:rational:approx:coeffs}
  &C(s,\kappa)   = \frac{2\kappa \sin(\pi s)}{\pi},\qquad
  a_l(s,\kappa) = \exp(2s(l-M_-(\kappa)-1)\kappa),\\[4pt]
  \label{eq:param:pde:coeffs}
  &b_l(s,\kappa) = \exp(2(l-M_-(\kappa)-1)\kappa),\qquad
  c_l(s,\kappa) = 1.
\end{align}
\end{subequations}
This specific choice of $\mathcal Q_s^{\kappa}$ satisfies \abrev{the error bound in~\eqref{eq:rational_scheme_convergence} with
\begin{equation}
  \varepsilon_s(\kappa) := \frac{2\sin(\pi s)}{\pi} \left[\frac{1}{2s} +
  \frac{1}{(2-2s)\lambda_0} \right]\left[\frac{\e^{-\pi^2/(4\kappa)}}{\sinh(\pi^2/(4\kappa))} + \e^{-\pi^2/(2\kappa)} \right]
  \label{eq:definition_epsilon}
\end{equation}}
which tends to zero exponentially fast as \abrev{$\kappa \to 0$; see~\cite[Remark 3.1]{2015_bonito}}.

We refer to~\cite{2017_bonitoa,2017_bonito,2021_bonito,2016_bonito} for applications of \abrev{the rational approximation given by~\eqref{eq:bp_rational_approximation}}
to the discretization of various types of fractional PDEs.

From now on, to simplify the notation, we will omit the dependence of \abrev{the summation limit $N$ and the coefficients}
$C,\, a_l,\, b_l,\, c_l$ in~\eqref{eq:definition_rational_function} on $s$ and $\kappa$.

\subsection{Finite element discretization}
\label{subsec:fe_discretization}

\abrev{As it was outlined above, in order to obtain a fully discrete approximation of the solution to problem~\eqref{eq:fractional_laplacian_strong}
one needs to discretize each parametric problem~\eqref{eq:continuous_parametric_problems}.
To that end, \abrev{similar to~\cite{2015_bonito, 2023_bulle},}
we use the FEM.
We note that the diffusion and/or reaction coefficients for parametric problems~\eqref{eq:continuous_parametric_problems} may
vary significantly between different problems;
for example, in the case of the BP rational approximation~\eqref{eq:bp_rational_approximation},
the diffusion coefficients $(b_l)_{l=1}^N$ vary from extremely small for  $l = 1$ to very large for $l = N$
 (cf.~\eqref{eq:param:pde:coeffs}).
Thus, as discussed in~\cite[Section 9.1.1]{2023_bulle},
using the same finite element mesh for all parametric problems may lead to wasted computational resources,
since some of these problems will be discretized on an over-refined mesh.
Therefore, the key idea and the main novelty of this study is to allow different meshes to be used
for the finite element discretization of different parametric problems in~\eqref{eq:continuous_parametric_problems}.
In particular, we will propose an algorithm that adaptively refines the finite element mesh \emph{individually} for each parametric problem and
employs an overlay of all the meshes (referred to as the \emph{union mesh}) in order to compute the
fully discrete approximation of the solution $u$ to the fractional Laplacian problem~\eqref{eq:fractional_laplacian_strong}.
}

For \abrev{$m \in \mathbb N_0$, let $\left\{\T_l^m,\ l = 1,\ldots,N\right\}$ be} a family of meshes \abrev{on $\Omega$,
where each mesh $\T_l^m$ is a conforming partition of $\Omega$ into compact non-degenerate simplices (cells).
The mesh $\T_l^m$ is associated with the $l$-th parametric problem in~\eqref{eq:continuous_parametric_problems},
and the subscript $m$ is the iteration counter in our adaptive algorithm, i.e.,
for each $l = 1,\ldots,N$, the mesh $\T^{m+1}_l$ is either
a refinement of the mesh $\T_l^m$ or $\T^{m+1}_l = \T_l^m$.}
%
\abrev{
For mesh refinement, we employ newest vertex bisection (NVB); see, e.g.,~\cite{stevenson,kpp}.
We assume that any mesh $\T_l^m$ employed for the discretization of parametric problems in~\eqref{eq:continuous_parametric_problems}
can be obtained by applying NVB refinement(s) to a given (coarse) initial mesh $\T^0$.
}

Given a mesh $\T_l^m$, we denote by $K$ a cell, by $F$
an edge (in \abrev{two dimensions}) or a face (in \abrev{three dimensions}),
and by $\chi$ a \abrev{vertex} of $\T_l^m$.
For a cell $K$, the set of edges/faces of $K$ is denoted by $\partial K$.
\abrev{Furthermore,} $\mathcal F_l^{m}$ \abrev{denotes} the set of \abrev{interior edges/faces} of $\T_l^m$,
\abrev{and for any edge/face $F \in \mathcal F_l^{m}$,
we denote by $n_F$ a unit normal vector to $F$.}
We assume that \abrev{these} normal vectors are fixed once and for all.
\abrev{Let $v$ be a sufficiently regular function.
For an edge/face $F \in \mathcal F_l^{m}$ shared by two cells $K_1$ and $K_2$,
we denote by $\llbracket v \rrbracket_F := v_{|K_1} - v_{|K_2}$ the jump of $v$ across $F$, and
let $\partial_{n_F} v := \nabla v \cdot n_F$ denote the normal derivative of $v$.}

Let $\omega \subseteq \Omega$, \abrev{$p \in \N$ and let} $\mathcal P^p(\omega)$ be
the space of polynomial functions of degree $p$ over $\omega$.
We denote by $\mathcal S^{m,p}_l$ the Lagrange finite element space of degree $p$ associated with the mesh $\T_l^m$:
\begin{equation} \label{eq:fe:space}
  \mathcal S^{m,p}_l := \left\{ v \in \abrev{H^1_0(\Omega)};\ v \in \mathcal P^p(K)\ \text{for all}\ K \in \T_l^m \right\}.
\end{equation}

For $l=1,\ldots,N$, we \abrev{introduce the Galerkin finite element formulation of the parametric problem~\eqref{eq:continuous_parametric_problems}:
find $w_l^m \in \mathcal S^{m,p}_l$ satisfying}
\begin{equation}
  b_l \scal{\nabla w_l^m}{\nabla v}_2 + c_l \scal{w_l^m}{v}_2 = \scal{f}{v}_2\quad \forall v \in \mathcal S^{m,p}_l.
  \label{eq:discrete_parametric_problems}
\end{equation}
\abrev{Thus, $w_l^m$ is the Galerkin} finite element approximation of the \abrev{solution $w_l$ to 
the parametric problem~\eqref{eq:continuous_parametric_problems}}.

\subsection{Union mesh}\label{subsec:union_mesh}

Let us consider a fixed \abrev{$m \in \N_0$ (i.e., a fixed iteration step of the adaptive algorithm that we are going to design)}.
\abrev{Recall that any two finite element meshes within the family $\left\{\T_l^m;\ l=1,\ldots,N \right\}$ might be different.}
In practice, in order to combine the \abrev{corresponding} finite element \abrev{approximations} $\left\{w_l^m;\ l=1,\ldots,N \right\}$,
we need to introduce a mesh $\widetilde{\T}^m$ such that the \abrev{associated} finite element space
$\widetilde{\mathcal S}^{m,p} := \big\{ v \in H^1_0(\Omega);\ v \in \mathcal P^p(K)\ \text{for all}\ K \in \widetilde\T^m \big\}$
contains all the functions $w_l^m$.
\abrev{Since each mesh in the family $\left\{\T_l^m;\ l=1,\ldots,N \right\}$ is either $\T^0$ or obtained by NVB refinements of $\T^0$,}
we can define the mesh $\widetilde{\T}^m$ as the coarsest common refinement \abrev{(i.e., the overlay)} of all meshes in \abrev{this family};
\abrev{we will call $\widetilde{\T}^m$ the \emph{union mesh}}.
Thus, \abrev{$w_l^m \in \mathcal S^{m,p} _l \subset \widetilde{\mathcal S}^{m,p}$} for all $l=1,\ldots,N$
and we define the \abrev{fully discrete approximation of the solution $u$ to problem~\eqref{eq:fractional_laplacian_strong}
as follows\ (cf.~\eqref{eq:rational_approximation_definition}):}
\begin{equation}
  \widetilde u_{\kappa}^m := C \sum_{l=1}^N a_l w_l^m \, \abrev{\in \widetilde{\mathcal S}^{m,p}}.
  \label{eq:discrete_approximation_u}
\end{equation}

\section{A posteriori error estimation} \label{sec:error:estimation}

The overall discretization error in the $L^2(\Omega)$-norm is given by
$\norm{u -  \widetilde u_{\kappa}^m}_2$,
where \abrev{$u$ is the solution to~\eqref{eq:fractional_laplacian_strong} and}
$\widetilde u_{\kappa}^m$ is defined in~\eqref{eq:discrete_approximation_u}.
Using the triangle inequality this error can be bounded as follows
\begin{equation}
  \norm{u - \widetilde u_{\kappa}^m}_2 \leqslant \norm{u - u_{\kappa}}_2 + \norm{u_{\kappa} - \widetilde u_{\kappa}^m}_2,
  \label{eq:fractional_error_decomposition}
\end{equation}
\abrev{where $u_\kappa$ is given by~\eqref{eq:rational_approximation_definition}}.
On the right-hand side of \eqref{eq:fractional_error_decomposition}, the contribution $\norm{u - u_{\kappa}}_2$
\abrev{to the overall discretization error} can be seen as the \emph{rational approximation error},
whereas the contribution $\norm{u_{\kappa} -  \widetilde u_{\kappa}^m}_2$ can be seen as the \emph{finite element discretization error}.
It turns out that, if $f \in L^2(\Omega)$ in \eqref{eq:fractional_laplacian_strong},
the rational approximation error \abrev{decays with}
the same rate as $\varepsilon_s(\kappa)$ in~\eqref{eq:rational_scheme_convergence}.
Indeed,
following the proof of~\cite[Theorem 3.5]{2015_bonito},
it is easy to show that
\begin{equation}
  \norm{u - u_{\kappa}}_2 \leqslant \varepsilon_s(\kappa) \norm{f}_2.
  \label{eq:convergence_rational_error}
\end{equation}

If the BP rational approximation~\eqref{eq:bp_rational_approximation}
is employed, then $\varepsilon_s(\kappa)$ is given by~\eqref{eq:definition_epsilon},
where the only unknown quantity is the lower bound $\lambda_0$ of the Laplacian spectrum.
Guaranteed lower bounds of the Laplacian spectrum can be computed, \abrev{see, e.g.,}~\cite{2017_cances,2021_carstensen}.
Therefore, in the case of the BP rational approximation,
the right-hand side of~\eqref{eq:convergence_rational_error} is fully computable
and provides an upper bound for the rational approximation error $\lVert u - u_{\kappa} \rVert_2$.
Thus, given a precision tolerance $\tol > 0$, we can choose the parameter $\kappa$
in~\eqref{eq:bp_rational_approximation} so that
\begin{equation*}
  \lVert u - u_{\kappa} \rVert_2  \leqslant \varepsilon_s(\kappa)\lVert f \rVert_2 \ll \tol.
\end{equation*}
In what follows we will assume that the rational approximation is sufficiently accurate 
so that the corresponding approximation error $\norm{u - u_{\kappa}}_2$ is negligible
compared to the finite element discretization error and
\begin{equation}
  \norm{u -  \widetilde u_{\kappa}^m}_2 \approx \norm{u_{\kappa} -  \widetilde u_{\kappa}^m}_2.
  \label{eq:neglect_rational_error} 
\end{equation}
Therefore, in this study, our focus is on the adaptive mesh refinement algorithm driven by
a posteriori estimates of the finite element error $\norm{u_{\kappa} -  \widetilde u_{\kappa}^m}_2$.
We refer to~\cite[Section 7.2]{2023_bulle} for details of a refinement algorithm that takes into account
both the finite element and the rational approximation~errors.

Our adaptive algorithm is steered by hierarchical a posteriori error estimators of the Bank--Weiser type~\cite{1985_bank}.
These estimators are computed using enriched finite element spaces on each cell of the mesh.

\subsection{Enriched local finite element spaces} \label{sec:enriched:fe:spaces}

For a cell $K \in \T_l^m$, we denote by $\mathcal S_l^{m,p}(K)$ the local Lagrange finite element space of degree $p$ over $K$.
The space $\mathcal S_l^{m,p}(K)$ is in fact the same as $\mathcal P^p(K)$.
Let $\widehat {\mathcal S}_l^m(K)$ 
denote an enriched local finite element space
that is obtained from ${\mathcal S}_l^{m,p}(K)$ by adding new basis functions.
Then, $\widehat {\mathcal S}_l^m(K)$ can be decomposed~as
\begin{equation}
  \widehat{\mathcal S}_l^m(K) = \mathcal S_l^{m,p}(K) \oplus \mathcal V_l^m(K),
  \label{eq:enrichment_space}
\end{equation}
where ${\mathcal V}_l^m(K) \subset H^1(K)$ and ${\mathcal V}_l^m(K) \cap {\mathcal S}_l^{m,p}(K) = \{0\}$.
The subspace ${\mathcal V}_l^m(K)$ in~\eqref{eq:enrichment_space}
is called the \emph{enrichment space}.
Similarly, for a cell $K$ of the union mesh $\widetilde \T^m$,
we denote the enrichment space by $\widetilde{\mathcal V}^m(K)$.

While enriched spaces 
can be constructed in many different ways
(see~\cite{1994_verfurth}), 
the most popular enrichment strategies are known as $p$-enrichment and $h$-enrichment.
In $p$-enrich\-ment, the space $\widehat{\mathcal S}_l^m(K)$ is defined as the finite element space
of a higher polynomial degree $\widehat{p} > p$, i.e.,
$
  \widehat{\mathcal S}_l^m(K) := \mathcal S_l^{m,\widehat p}(K).
$
In this case, the enrichment space $\mathcal V_l^m(K)$ in~\eqref{eq:enrichment_space} is
the space of \emph{polynomials} of degree $\widehat p$ defined over $K$ and
vanishing at the degrees of freedom of $\mathcal S_l^{m,p}(K)$.
Specific examples of $p$-enrichment can be found in~\cite{1985_bank}.

In $h$-enrichment, the space $\widehat{\mathcal S}_l^m(K)$ is defined
as the Lagrange finite element space of degree $p$
associated with a uniform partition of the cell $K$.
In this case, the space $\mathcal V_l^m(K)$ is the space of (continuous) \emph{piecewise polynomials}
of degree $p$ defined over a uniform partition of $K$ and
vanishing at the degrees of freedom of $\mathcal S_l^{m,p}(K)$.
We refer to~\cite{2000_ainsworth} for specific examples of $h$-enrichment.

\subsection{Local error estimation} \label{sec:local_estimation}

Let $K$ be a cell of $\T_l^m$ for some $l=1,\ldots,N$.
In the Bank--Weiser error estimation strategy, the enrichment space ${\mathcal V}_l^m(K)$ in~\eqref{eq:enrichment_space}
is utilized to obtain a local estimator approximating the finite element error on the cell~$K$.
More precisely, this error estimator is defined as the solution to a local Neumann problem on $K$.
The right-hand side of this problem is written it terms of the interior and edge/face residuals
associated with the Galerkin approximation $w_l^m$ satisfying~\eqref{eq:discrete_parametric_problems}.
The interior and edge/face residuals are denoted by $r_{l,K}^m$ and $J_{l,F}^m$, respectively, and defined as follows:
\begin{equation}
  r_{l,K}^m := (f - c_l w_l^m + b_l \Delta w_l^m)|_{K}
  \label{eq:interior_residual}
\end{equation}
and
\begin{equation}
  J_{l,F}^m :=
    \begin{cases}
      b_l \llbracket \partial_{n_F} w_l^m \rrbracket_F & \text{if } F \in \mathcal F_l^{m},\\
      0                                          & \text{otherwise.}
    \end{cases}
    \label{eq:edges_residual}
\end{equation}
Then, the local Neumann problem associated with $K$ reads as:
find $e_{l,K}^m \in \mathcal V_l^m(K)$ such that for all $v \in \mathcal V_l^m(K)$ there holds
\begin{equation}
  b_l \scal{\nabla e_{l,K}^m}{\nabla v}_{2,K} + c_l \scal{e_{l,K}^m}{v}_{2,K} =
  \scal{r_{l,K}^m}{v}_{2,K} - \frac{1}{2} \sum_{F \in \partial K} \scal{J_{l,F}^m}{v}_{2,F}.
  \label{eq:local_parametric_bw_equation}
\end{equation}
The function $e_{l,K}^m$ is an approximation of the local error $(w_l-w_l^m)|_{K}$ and, consequently,
we consider the following \rev{heuristic} \emph{local error indicators}:
\begin{equation}
  \eta_{l,K}^m := \norm{e_{l,K}^m}_{2,K} \approx \norm{w_l-w_l^m}_{2,K},
  \label{eq:definition_local_estimator}
\end{equation}
will be used in our adaptive algorithm to mark the cells of $\mathcal T_l^m$ for refinement.

\subsection{Global error estimation} \label{sec:global:estimation}

In the previous section we have introduced the \emph{local}
error indicators that will be guiding adaptive refinement of finite element meshes.
Our goal now is to obtain a computable estimate 
of the \emph{global} finite element error $\norm{u_{\kappa} -  \widetilde u_{\kappa}^m}_2$,
where $u_{\kappa}$ is given by~\eqref{eq:rational_approximation_definition}
and $\widetilde u_{\kappa}^m$ is defined in~\eqref{eq:discrete_approximation_u}.
Global error estimates are used to control the overall (finite element) error across the computational domain and,
for a given tolerance, they provide a stopping criterion in adaptive algorithms.

If the same finite element mesh is employed for all parametric problems~\eqref{eq:discrete_parametric_problems},
i.e., if $\T_l^m = \T^m$ for all $l=1,\ldots,N$, then obtaining a global error estimate is straightforward.
Indeed, for each cell $K \in \T^m = \widetilde \T^m$,
the error estimators $e_{l,K}^m$ defined by~\eqref{eq:local_parametric_bw_equation}
can be combined across all parametric problems and one can define
\begin{equation}
  \widetilde e_{K}^m := C \sum_{l=1}^N a_l e_{l, K}^m
  \label{eq:local_fractional_bw_solution1}
\end{equation}
as an approximation of the local error $(u_\kappa - \widetilde u_{\kappa}^m)|_{K}$; 
cf.~\eqref{eq:rational_approximation_definition},~\eqref{eq:discrete_approximation_u}.
Hence, the global error estimate $\widetilde \eta^m$ can be easily defined as follows:
\begin{equation}
   \norm{u_{\kappa}-\widetilde u_{\kappa}^m}_2^2 =
   \sum_{K \in {\mathcal T}^m} \norm{u_{\kappa}-\widetilde u_{\kappa}^m}_{2,K}^2
   \approx 
   \sum_{K \in {\mathcal T}^m} \norm{\widetilde e_{K}^m}_{2,K}^2
   =: (\widetilde \eta^m)^2.
   \label{eq:union_mesh_fractional_estimator1}
\end{equation}
This approach to the error estimation within the \emph{single-mesh} discretization framework for the spectral fractional Laplacian
has been first proposed in~\cite{2023_bulle}.
The numerical experiments included in~\cite{2023_bulle} have demonstrated the effectivity of this approach.

In our \emph{multimesh} discretization framework,
there is no straightforward way to combine local error estimators $e_{l,K}^m$ over parametric problems,
simply because the meshes underlying finite element approximations for different
parametric problems may not share a particular cell $K$.
Notably, the definition of $\widetilde e_{K}^m$ in~\eqref{eq:local_fractional_bw_solution1} is no longer valid in this setting.
We propose two approaches to address this issue in deriving a computable global error estimate
in the multimesh setting.

\subsubsection{Global error estimation based on the triangle inequality}
\label{sec:global:estimate:1}
The first approach relies on the triangle inequality to obtain a bound on the true error as follows:
\begin{equation}
  \norm{u_{\kappa} - \widetilde u_{\kappa}^m}_2
  = \left \lVert C \sum_{l=1}^N a_l w_l - C \sum_{l=1}^N a_l w_l^m\right \rVert_2 
  \leqslant C \sum_{l=1}^N a_l \norm{w_l - w_l^m}_2.
  \label{eq:bound_error_triangle_inequality}
\end{equation}
Using the local error estimates in~\eqref{eq:definition_local_estimator},
we define the global error estimates $\{\eta_l^m;\ l=1,\ldots,N\}$ independently
for each parametric problem as follows:
\begin{equation}
  (\eta_l^m)^2 := \sum_{K \in \mathcal T_l^m} (\eta_{l,K}^m)^2\quad \forall l=1,\ldots,N.
  \label{eq:global_parametric_estimator}
\end{equation}
Thus, from~\eqref{eq:definition_local_estimator} we have
\begin{equation}
  \norm{w_l - w_l^m}_2 \approx \eta_l^m.
  \label{eq:local_error_estimator}
\end{equation}
Then the overall global error estimate is defined as
\begin{equation}
  \eta^m := C \sum_{l=1}^N a_l \eta_l^m
  \overset{\eqref{eq:global_parametric_estimator}}{=}
  C \sum_{l=1}^N a_l \bigg[\sum_{K \in \mathcal T_l^m} (\eta_{l,K}^m)^2 \bigg]^{1/2}
  \label{eq:triangular_inequality_fractional_estimator}
\end{equation}
and there holds
\begin{equation}
  \norm{u-\widetilde u_{\kappa}^m}_2
  \overset{\eqref{eq:neglect_rational_error}}{\approx}
  \norm{u_{\kappa} - \widetilde u_{\kappa}^m}_2
  \overset{\eqref{eq:bound_error_triangle_inequality}}{\leqslant}
  C \sum_{l=1}^N a_l \norm{w_l - w_l^m}_2 
  \overset{\eqref{eq:local_error_estimator}}{\approx}
  C \sum_{l=1}^N a_l \eta_l^m = \eta^m.
  \label{eq:true_error_bound_triangular_inequality_estimator}
\end{equation}
This approach has the advantage of being very easy to implement.
While $\eta^m$ given by~\eqref{eq:triangular_inequality_fractional_estimator} is cheap to compute from the local
error indicators $\{\eta_{l,K}^m;\ K \in \mathcal T_l^m,\ l=1,\ldots,N\}$,
the use of the triangle inequality in~\eqref{eq:bound_error_triangle_inequality}
may affect the effectivity of the resulting global error estimate.
We will address this issue when we discuss numerical results in~\S\ref{sec:numerical_results}.

\subsubsection{Global error estimation based on the union mesh} \label{sec:global:estimate:2}

In our second approach,
instead of using the error indicators $\eta_{l,K}^m$ from the meshes $\T_l^m$ ($l=1,\ldots,N$),
we compute the local error estimators directly on each cell $\widetilde K$ of
the \emph{union mesh} $\widetilde \T^m$.
In this case, the local Neumann problem reads as follows: find
$e_{l,\widetilde K}^m \in \widetilde{\mathcal V}^m(\widetilde K)$ such that
for all $v \in \widetilde{\mathcal V}^m(\widetilde K)$ there holds
\begin{equation} \label{eq:local_parametric_bw_equation:union_mesh}
  b_l \! \scal{\nabla e_{l,\widetilde K}^m}{\nabla v}_{2,\widetilde K} +
  c_l \! \scal{e_{l,\widetilde K}^m}{v}_{2,\widetilde K} =
  \! \scal{r_{l,\widetilde K}^m}{v}_{2,\widetilde K} -
  \frac{1}{2} \sum_{\widetilde F \in \partial \widetilde K} \! \scal{J_{l,\widetilde F}^m}{v}_{2,\widetilde F},
\end{equation}
where, for each $l=1,\ldots,N$, the residuals $r_{l,\widetilde K}^m$ and $J_{l,\widetilde F}^m$
(for each $\widetilde F \in \partial \widetilde K$)
are computed similarly to \eqref{eq:interior_residual} and \eqref{eq:edges_residual} from the parametric solutions $w_l^m$
interpolated at the degrees of freedom associated with the union mesh $\widetilde \T^m$.
Thus, for each cell $\widetilde K$ of the union mesh, we obtain an approximation
$e_{l,\widetilde K}^m \in \widetilde{\mathcal V}^m(\widetilde K)$ of the local error $(w_l - w_l^m)|_{\widetilde K}$.
Then, in analogy with~\eqref{eq:local_fractional_bw_solution1}, we define
\begin{equation*}
  \widetilde e_{\widetilde K}^m := C \sum_{l=1}^N a_l \widetilde e_{l,\widetilde K}^m
\end{equation*}
as an approximation of the local error
$(u_\kappa - \widetilde u_{\kappa}^m)|_{\widetilde K}$ for each $\widetilde K \in \widetilde \T^m$.
Then
$\norm{u_{\kappa} - \widetilde u_{\kappa}^m}_{2,\widetilde K} \approx
  \norm{\widetilde e_{\widetilde K}^m}_{2,\widetilde K}
$
and defining the overall global error estimate as
\begin{equation}
  \widetilde \eta^m :=
  \Bigg[
  \sum_{\widetilde K \in \widetilde{\mathcal T}^m} \norm{\widetilde e_{\widetilde K}^m}_{2,\widetilde K}^2
  \Bigg]^{1/2},
  \label{eq:union_mesh_fractional_estimate_def}
\end{equation}
there holds (cf.~\eqref{eq:union_mesh_fractional_estimator1} in the single-mesh case)
\begin{equation*}
  \norm{u-\widetilde u_{\kappa}^m}_2
  \overset{\eqref{eq:neglect_rational_error}}{\approx}
  \norm{u_{\kappa} - \widetilde u_{\kappa}^m}_2
  =
  \Bigg[\sum_{\widetilde K \in \widetilde{\mathcal T}^m}
  \norm{u_{\kappa} - \widetilde u_{\kappa}^m}_{2,\widetilde K}^2 \Bigg]^{1/2}
  \approx \widetilde \eta^m.
\end{equation*}
This approach avoids the triangle inequality in \eqref{eq:bound_error_triangle_inequality}
at the expense of being more computationally demanding.
It requires the construction of the union mesh $\widetilde{\mathcal T}^m$,
the interpolation of each parametric solution $w_l^m$
at the degrees of freedom associated with $\widetilde{\mathcal T}^m$, 
and the solution of local Neumann problems~\eqref{eq:local_parametric_bw_equation:union_mesh}
on the cells of $\widetilde{\mathcal T}^m$.
On the other hand, in the adaptive algorithm, the global error estimate only needs to be computed periodically
in order to check whether the stopping criterion is satisfied.
Therefore,
in order to reduce the overall computational complexity of the algorithm, we
compute the global error estimate $\widetilde\eta^m$ given by~\eqref{eq:union_mesh_fractional_estimate_def}
and check the stopping criterion only at iterations $m \in \{0, k, 2k, 3k, \ldots\}$,
where $k \in \mathbb N$ is fixed.

\section{Adaptive multimesh refinement algorithm}\label{sec:adaptive:refinement}

The adaptive refinement algorithm
generates a sequence of fully discrete approximations
$\{\widetilde u_{\kappa}^m;\ m\in \N_0\}$ by iterating the following loop:
\begin{equation}
  \textbf{Solve} \Longrightarrow \textbf{Estimate} \Longrightarrow \textbf{Mark} \Longrightarrow \textbf{Refine}.
  \label{eq:adaptive_refinement_loop}
\end{equation}
The ingredients of the modules \textbf{Solve} and \textbf{Estimate}
were described in~\S\ref{sec:multimesh:scheme} and \S\ref{sec:error:estimation}, respectively.
While these two modules have to be executed for each
finite element formulation in~\eqref{eq:discrete_parametric_problems}
on the corresponding underlying mesh,
the modules \textbf{Mark} and \textbf{Refine}, discussed below,
are designed to optimize the refinement \emph{across all the meshes}.

\subsection{Multimesh marking}\label{subsec:marking}

The core idea underpinning our marking strategy in the multimesh setting
is to consider the local error indicators $\eta_{l,K}^m$ (see~\eqref{eq:definition_local_estimator})
from all parametric problems and apply a marking algorithm to the joint set of \emph{weighted} error indicators
$\big\{a_l \eta_{l,K}^m;\ K \in \mathcal T_l^m,\ l=1,\ldots,N\big\}$,
where $a_l$ are the coefficients in the rational approximation~\eqref{eq:definition_rational_function}
(see also the representations of the semi-discrete and fully discrete approximations
in~\eqref{eq:rational_approximation_definition} and~\eqref{eq:discrete_approximation_u}, respectively).

In this study we use D{\"o}rfler marking~\cite{1996_dorfler}, but
other marking strategies could also be applied in a similar way.
Let $\theta \in (0,1]$ be a fixed marking threshold.
\rev{The module~\textbf{Mark} in the adaptive loop~\eqref{eq:adaptive_refinement_loop}
finds a family of sets $\{\mathcal M_l^m \subseteq \T_l^m;\ l=1,\ldots,N\}$} such that
\begin{equation}
  \theta\,
  \sum_{l=1}^N \sum_{K \in \mathcal T_l^m} \big(a_l \eta_{l,K}^m\big)^2
  \leqslant \sum_{l=1}^N \sum_{K \in \mathcal M_l^m} \big(a_l \eta_{l,K}^m\big)^2,
  \label{eq:marking}
\end{equation}
with the cumulative cardinality $\sum_{l=1}^N \#\mathcal M_l^m$ that is minimized over all the sets
that satisfy~\eqref{eq:marking}.

\subsection{Multimesh refinement} \label{subsec:refinement}

The marking sets $\mathcal M^m_l$ generated by the module~\textbf{Mark}
are fed into the \textbf{Refine} module that
performs local NVB refinement of individual marked cells in each finite element mesh $\T_l^m$, $l = 1,\ldots,N$
(note that some cells which have not been marked but are adjacent to the marked cells
will also need to be bisected in order to ensure the conformity of the resulting refined mesh).

We emphasize here that coefficients $a_l$ (see, e.g.,~\eqref{eq:rational:approx:coeffs})
play an important role in the proposed marking and refinement strategy.
Used as weights in D{\"o}rfler marking (see~\eqref{eq:marking}),
they amplify or diminish contributions of
the local error indicators associated with the mesh $\T_l^m$ according to the significance
of the corresponding $l$-th term in the rational approximation~\eqref{eq:definition_rational_function}
and in representation~\eqref{eq:discrete_approximation_u}.
For example, if a coefficient $a_l$ is very small,
then the $a_l$-weighted contributions of the error indicators $\eta_{l,K}^m$ ($K \in \T_l^m$)
in~\eqref{eq:marking} will be insignificant, and
it is very likely that the corresponding marking set $\mathcal M_l^m$ will be empty.
This means that no cell from the mesh $\mathcal T_l^m$ will be selected for refinement.
Hence, this mesh $\mathcal T_l^m$
as well as the corresponding finite element approximation and
the associated local error indicators will all carry over to the next iteration of the adaptive loop,
i.e.,~$\T_l^{m+1} := \T_l^m$,
$w_l^{m+1} := w_l^m$, and
$\eta_{l,K}^{m+1} := \eta_{l,K}^m$ for all $K \in \mathcal T_l^{m+1} = \mathcal T_l^m$.
Therefore, one will not need to run the solve or error estimation routines
for the $l$-th parametric problem during the $(m+1)$-st iteration.
Thus, even though rational approximation schemes involve
solving many parametric problems~\eqref{eq:discrete_parametric_problems},
employing the multimesh framework allows to develop an adaptive algorithm
where only a small number of these problems are actually solved at each iteration
of the adaptive loop
(except the very first iteration, where all parametric problems have to be solved
on the coarsest mesh).
We will illustrate this aspect of the multimesh approach when we present
numerical results in~\S\ref{sec:numerical_results}.

\subsection{Adaptive algorithm} \label{sec:algorithm}

In this section, we present a generic adaptive algorithm
for computing multimesh rational approximations of the solution to
problem~\eqref{eq:fractional_laplacian_strong}.

The algorithm takes as inputs
the fractional power $s \in (0,1)$,
an initial (coarse) mesh $\T^0$, 
the D{\"o}rfler marking parameter $\theta \in (0,1]$,
the stopping tolerance $\tol$,
and the counter $k \in \N$ that determines the iterations at which the stopping criterion is checked
(these are iterations $0, k, 2k, 3k,\ldots$).
For some $m_* = j k$ with $j \in \N_0$, the algorithm generates a sequence
$\big\{ \widetilde u_{\kappa}^m;\ m=0, k, 2k,\ldots,m_* \big\}$ of fully discrete approximations
(computed from parametric Galerkin approximations $w_l^m \in \mathcal S^{m,p}_l$
using~\eqref{eq:discrete_approximation_u}) and two sequences of global a posteriori error estimates,
$\big\{ \eta^m;\ m=0, k, 2k,\ldots,m_* \big\}$ and $\big\{ \widetilde \eta^m;\ m=0, k, 2k,\ldots,m_* \big\}$,
where $\eta^m$ and $\widetilde\eta^m$ are defined by~\eqref{eq:triangular_inequality_fractional_estimator}
and~\eqref{eq:union_mesh_fractional_estimate_def}, respectively.

The adaptive algorithm is listed in~\Cref{alg:algorithm_outline}.
It contains the following eight subroutines:

\begin{itemize}
\item \texttt{GenerateRationalScheme}($s,\tol$)---a subroutine that
computes the values of the coefficients $C$, $a_l$, $b_l$, $c_l$ ($l=1,\ldots,N$)
of the rational approximation $\mathcal Q^{\kappa}_s$ in~\eqref{eq:definition_rational_function}
for a given fractional power $s$ of the Laplacian and a prescribed tolerance $\tol$.
As discussed in~\S\ref{sec:error:estimation},
the parameter $\kappa$ in the rational approximation~\eqref{eq:definition_rational_function} is chosen such that
the error bound $\varepsilon_s(\kappa)$ in~\eqref{eq:rational_scheme_convergence}
is negligible compared to~$\tol$.
%

\item \texttt{Solve}\big($\T^m_l,b_l,c_l$\big)---a subroutine that generates the Galerkin approximation
$w_l^m \in \mathcal S^{m,p}_l$ satisfying~\eqref{eq:discrete_parametric_problems}.

\item \texttt{Estimate}\big($w^m_l,\T^m_l,b_l,c_l$\big)---a subroutine that
computes local error indicators $\eta_{l,K}^m$ ($K \in \T^m_l$) in~\eqref{eq:definition_local_estimator} by
solving the local Neumann problem~\eqref{eq:local_parametric_bw_equation};
see the discussion in~\S\ref{sec:local_estimation}.

\item \texttt{Mark}\big($\{ \T_l^m, \{\eta_{l,K}^m\}_{K \in \T_l^m}, a_l \}_{l=1}^N, \theta$\big)---a subroutine
that generates the marking sets $\mathcal M_l^m$ for all finite element meshes $\T_l^m$, $l = 1,\ldots,N$.
The sets $\mathcal M_l^m$ are generated using the D{\"o}rfler marking strategy adapted to the multimesh framework;
see~\eqref{eq:marking}.

\item \texttt{Refine}\big($\T_l^m, \mathcal M^m_l$\big)---a subroutine that generates a refined mesh $\T_l^{m+1}$
as described in~\S\ref{subsec:refinement}.

\item \texttt{GenerateUnionMesh}\big($\{\T^m_l\}_{l=1}^N$\big)---a subroutine that generates
the union mesh $\widetilde \T^m$; 
see~\S\ref{subsec:union_mesh}.

\item \texttt{GlobalEstimate}\big($C,\! \{ w_l^m, \{\eta_{l,K}^m\}_{K \in \T_l^m}, a_l, b_l, c_l \}_{l=1}^N,\!
                                                      \widetilde \T^m$\big)---a subroutine
that computes two global error estimates $\eta^m$ and $\widetilde \eta^m$
defined in~\eqref{eq:triangular_inequality_fractional_estimator}
and~\eqref{eq:union_mesh_fractional_estimate_def}, respectively;
see~\S\ref{sec:global:estimation} for a detailed discussion.
In the algorithm, the error estimate $\widetilde \eta^m$ is used in the stopping criterion on line~14.

\item \texttt{SolutionOnUnionMesh}\big($C,\! \{w^m_l, a_l\}_{l=1}^N,\! \widetilde \T^m$\big)---a subroutine that
computes the fully discrete approximation $\widetilde u_{\kappa}^m$
from the Galerkin approximations $\{w^m_l\}_{l=1}^N$ interpolated at the nodes
(and, if $p > 1$, at other degrees of freedom) on the union mesh $\widetilde \T^m$; see~\eqref{eq:discrete_approximation_u}.
\end{itemize}

\begin{algorithm}[!th]
  \caption{Multimesh finite element discretization and adaptive refinement}
  \begin{algorithmic}[1]
    \Require $s \in (0,1)$, $\T^0$, $\theta \in (0,1]$, $\tol > 0$, $k \in \N$
    \State $\big\{ C,\ \{a_l,b_l,c_l\}_{l=1}^N \big\} \leftarrow$ \texttt{GenerateRationalScheme}($s,\tol$)
    \State $\T^0_l \leftarrow \T^0\ \forall l =1,\ldots,N$
    \For{$m=0,1,\ldots$}
      \For{$l=1,\ldots,N$}
        \If{$m=0$ \textbf{or} ($m \ge 1$ \textbf{and} $\mathcal T^{m-1}_l \neq \mathcal T^m_l$)}
          \State $w^m_l$ $\leftarrow$ \texttt{Solve}\big($\T^m_l,b_l,c_l$\big)
          \State $\{\eta_{l,K}^m\}_{K \in \mathcal T^m_l}$ $\leftarrow$
          \texttt{Estimate}\big($w^m_l,\T^m_l,b_l,c_l$\big)
        \EndIf
      \EndFor
      \If{$m \!\mod k = 0$}
        \State $\widetilde \T^m$ $\leftarrow$ \texttt{GenerateUnionMesh}\big($\{\T^m_l\}_{l=1}^N$\big)
        \State $\eta^m$, $\widetilde \eta^m$ $\leftarrow$
        \texttt{GlobalEstimate}\big($C, \{ w_l^m, \{\eta_{l,K}^m\}_{K \in \T_l^m}, a_l, b_l, c_l \}_{l=1}^N,
                                                     \widetilde \T^m$\big)
        \State $\widetilde u_{\kappa}^m$ $\leftarrow$
        \texttt{SolutionOnUnionMesh}\big($C, \{w^m_l, a_l\}_{l=1}^N, \widetilde \T^m$\big)
        \If{$\widetilde \eta^m < \tol$}
          \State \textbf{break} 
        \EndIf
      \EndIf
      \State $\{\mathcal M^m_l\}_{l=1}^N$ $\leftarrow$
      \texttt{Mark}\big($\{ \T_l^m, \{\eta_{l,K}^m\}_{K \in \T_l^m}, a_l \}_{l=1}^N, \theta$\big)
      \For{$l =1,\ldots,N$}
        \If{$\mathcal M^m_l = \varnothing$}
          \State $\T^{m+1}_l \leftarrow \T^m_l$
          \State $w^{m+1}_l \leftarrow w_l^m$
          \For{$K \in \T^{m+1}_l$}
            \State $\eta^{m+1}_{l,K} \leftarrow \eta_{l,K}^m$
          \EndFor
        \Else
          \State $\T^{m+1}_l$ $\leftarrow$ \texttt{Refine}\big($\T_l^m, \mathcal M^m_l$\big)
        \EndIf
      \EndFor
    \EndFor
    \Ensure $\big\{ \widetilde u_{\kappa}^{m},\, \eta^{m},\, \widetilde \eta^{m};\ m=0, k, 2k,\ldots,m_* \big\}$
    for some $m_* = j k$ with $j \in \N_0$
  \end{algorithmic}
  \label{alg:algorithm_outline}
\end{algorithm}

We recall from the discussion in~\S\ref{subsec:refinement} that
not all finite element meshes will be refined at each iteration of the adaptive algorithm.
Indeed, if for a given $l \in \{1,\ldots,N\}$, the marking set $\mathcal M^m_l$ is empty,
then the mesh $\T^m_l$ is not refined at this iteration (i.e.,~$\T_l^{m+1} := \T_l^m$), and
all the quantities associated with $\T^m_l$ are carried over to the next iteration
(see lines~20--25 in Algorithm~\ref{alg:algorithm_outline})
saving significant computational resources and reducing the computational time.

Since the parametric problems are independent from each other,
we emphasize that the \textbf{for}-loops in $l$ (see lines~4--9 and~19--29 in Algorithm~\ref{alg:algorithm_outline})
can be parallelized, meaning that each of the subroutines \texttt{Solve}, \texttt{Estimate}, and \texttt{Refine}
is run simultaneously for all parametric problems.
In addition, given any pair of cells $K$ and $K'$ in $\T_l^m$, the corresponding local
problems~\eqref{eq:local_parametric_bw_equation} for $K$ and $K'$ are also independent from each other.
Therefore, the computation of local a posteriori error indicators in the subroutine \texttt{Estimate}
can be vectorized over cells, further reducing the computational time.

\section{Numerical results} \label{sec:numerical_results}

The aim of this section is threefold.
Firstly, we will investigate the effectivity and robustness
of two strategies for global error estimation in the multimesh setting
(see~\S\ref{sec:global:estimation}).
Secondly, we will demonstrate the performance of the adaptive algorithm
presented in~\S\ref{sec:algorithm}.
Finally, we will compare the performance of adaptive algorithms (in terms of convergence rates,
computational complexity, and runtime) in the single-mesh and multimesh settings.
To that end, we will consider three representative test problems
for the spectral fractional Laplacian in two dimensions.

The numerical results presented here were produced using a MATLAB implementation of
Algorithm~\ref{alg:algorithm_outline} within the finite element toolbox
T-IFISS~\cite{2021_tifiss_paper, tifiss_software}.
The following implementation details are worth noting:

\begin{itemize}

\item
in the subroutine \texttt{GenerateRationalScheme},
the Bonito--Pasciak rational approximation
(see~\eqref{eq:bp_rational_approximation}--\eqref{eq:definition_epsilon})
is used;
in all computations, we set $\kappa = 0.26$
in~\eqref{eq:bp_rational_approximation}
so that the error bound $\varepsilon_s(\kappa)$ in~\eqref{eq:definition_epsilon}
does not exceed $2 \cdot 10^{-8}$ for all test cases;
this guarantees that the rational approximation error in~\eqref{eq:convergence_rational_error}
is significantly smaller than the finite element discretization error and~\eqref{eq:neglect_rational_error} holds;

\item
the first-order (P1) finite element approximations are employed (i.e., $p=1$ in~\eqref{eq:fe:space});

\item
the implementation of the a posteriori error estimation strategy in~\S\ref{sec:error:estimation}
employs the enrichment space $\mathcal V_l^m(K)$ in~\eqref{eq:enrichment_space}
that is obtained via $h$-enrichment;
more precisely, $\mathcal V_l^m(K)$ contains piecewise linear functions over subelements
obtained by uniform (NVB) refinement of $K$
(see, e.g.,~\cite[Figure~5.2]{2000_ainsworth}
and~\cite[Section~2.1]{2021_tifiss_paper});

\item
\rev{the union mesh $\widetilde \T^m$ is assembled via the Delaunay triangulation method with edge constraints; 
specifically, the nodes and edges of all meshes in the family $\{\mathcal T_l^m;\ l=1,\ldots,N\}$
are passed to a Delaunay triangulation algorithm that generates the overlay of the meshes in this family;
}

\item
\rev{we follow~\cite{1996_dorfler} and implement the marking strategy in~\eqref{eq:marking}
by sorting the set
$\{a_l \eta_{l,K}^m;$ $K \,{\in}\, \mathcal T_l^m,\ l \,{=}\, 1,\ldots,N\}$ in decreasing order and
iteratively selecting the largest (weighted) error indicators until~\eqref{eq:marking} is reached. 
This strategy is of log-linear complexity with respect to the cardinality of the above set. 
More efficient strategies that do not rely on sorting and are of linear complexity
can also be employed here, see, e.g.,~\cite{2020_pfeiler} and the references therein.}

\item
local mesh refinement is performed using a variant of the newest vertex bisection method
called the \emph{longest edge bisection}; see~\cite[Section~2.1]{2021_tifiss_paper} and the references therein for details.

\end{itemize}

Unless specified otherwise, in the numerical experiments presented below, we set
the marking threshold (see~\eqref{eq:marking}) to $\theta = 0.5$
and set $k=1$ (i.e., the stopping criterion for the global error estimate $\widetilde\eta^m$
is checked at every iteration of the adaptive loop).

\subsection{Test case I: square domain, unit right-hand side function}\label{sec:test:case:I}

Let us set $f  \equiv 1$ and consider the model problem~\eqref{eq:fractional_laplacian_strong}
on the square domain $\Omega = (-1,1)^2$.
The primary aim of this test case is to investigate the effectivity and robustness of two global error estimation strategies
introduced in~\S\ref{sec:global:estimation} in the multimesh setting.
To that end, for $s \in \{0.3,\, 0.5,\, 0.7\}$, we set $\T^0$ as a uniform mesh of 512 right-angled triangles
and run our adaptive multimesh refinement algorithm.
In addition, for each $s$ from the same range, we generate a reference (fully discrete) solution
$u_{\text{ref}} \in {\mathcal S}^{\text{ref}} \subset H^1_0(\Omega)$
to problem~\eqref{eq:fractional_laplacian_strong}
using a single, highly refined Shishkin-type mesh $\T^{\text{ref}}$.
Such meshes contain anisotropic elements in the boundary layer,
they are commonly used for the numerical solution of singularly perturbed differential equations
and provide accurate approximations (see, e.g.,~\cite{KoptevaR_10_SMN} and the references therein). 
We use the reference solutions to compute the effectivity indices
for global a posteriori error estimates
$\eta^m$ and $\widetilde\eta^m$ defined by~\eqref{eq:triangular_inequality_fractional_estimator}
and~\eqref{eq:union_mesh_fractional_estimate_def}, respectively.
The effectivity indices are defined as follows:
\begin{equation}
   \Theta^m := \frac{\eta^m}{\norm{u_{\text{ref}}-\widetilde u_{\kappa}^m}_2},
   \quad
   \widetilde\Theta^m := \frac{\widetilde\eta^m}{\norm{u_{\text{ref}}-\widetilde u_{\kappa}^m}_2},
   \label{eq:effectivity:indices}
\end{equation}
where $\widetilde u_{\kappa}^m$ is the fully discrete approximation generated at iteration $m$
of Algorithm~\ref{alg:algorithm_outline}.

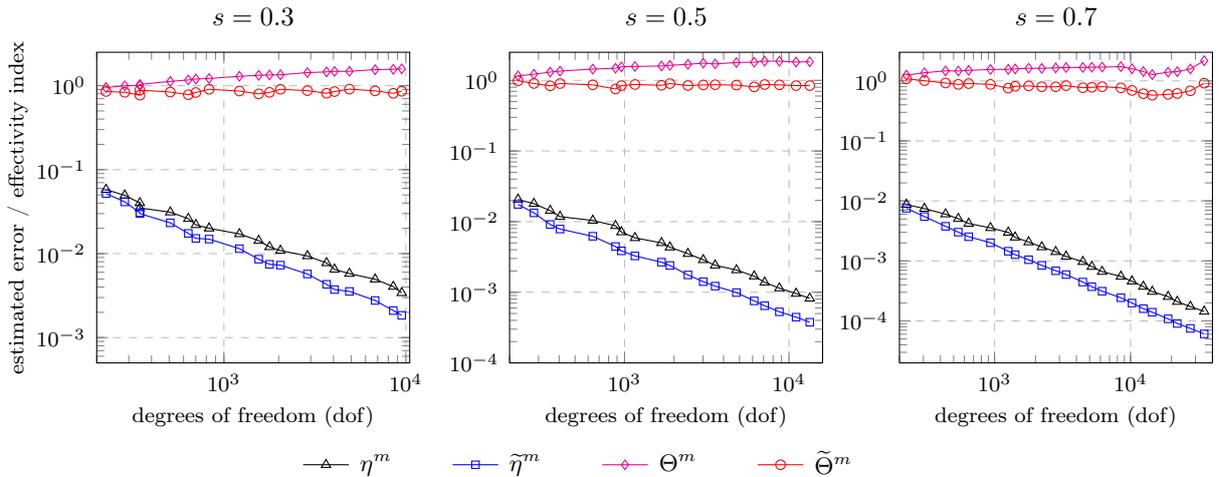
\begin{figure}[t!]
\begin{tikzpicture}
\pgfplotstableread{tp2_multimesh_s03.dat}{\one}
\begin{loglogaxis}
[
title={\footnotesize $s = 0.3$},
width = 4.8cm, height = 4.8cm,						
xlabel={degrees of freedom (dof)}, 					
xlabel style={font=\fontsize{8pt}{10pt}\selectfont},		
ylabel={estimated error / effectivity index},				
ylabel style={font=\fontsize{8pt}{10pt}\selectfont}, 		
ymajorgrids=true, xmajorgrids=true, grid style=dashed,	
xmin=200, xmax=1.05*10^4,						
ymin = 5*10^(-4),	 ymax = 2.5,							
legend columns=4,
legend entries={
{$\eta^m$\qquad\qquad},
{$\widetilde\eta^m$\qquad\qquad},
{$\Theta^m$\qquad\qquad},
{$\widetilde\Theta^m$\qquad\qquad}
},
legend to name=legend_tp2_1,
legend style={legend cell align=left, row sep=5pt, column sep=2pt, fill=none, draw=none, font={\fontsize{9pt}{12pt}\selectfont}},
]
\addplot[black,mark=triangle,mark size=2pt]		table[x=intdofs, y=comperror3]{\one};
\addplot[blue,mark=square,mark size=1.4pt]		table[x=intdofs, y=comperror2_union]{\one};
\addplot[magenta,mark=diamond,mark size=1.8pt]		table[x=intdofs, y=eff_index3]{\one};
\addplot[red,mark=o,mark size=1.8pt]		table[x=intdofs, y=eff_index2_union]{\one};
%
\end{loglogaxis}
\end{tikzpicture}
\hspace{-5pt}
\begin{tikzpicture}
\pgfplotstableread{tp2_multimesh_s05.dat}{\two}
\begin{loglogaxis}
[
title={\footnotesize $s = 0.5$},
width = 4.8cm, height = 4.8cm,						
xlabel={degrees of freedom (dof)}, 					
xlabel style={font=\fontsize{8pt}{10pt}\selectfont},		
ymajorgrids=true, xmajorgrids=true, grid style=dashed,	
xmin=200, xmax=16000,						
ymin = 1*10^(-4),	 ymax = 2.5,							
]
\addplot[black,mark=triangle,mark size=2pt]		table[x=intdofs, y=comperror3]{\two};
\addplot[blue,mark=square,mark size=1.4pt]		table[x=intdofs, y=comperror2_union]{\two};
\addplot[magenta,mark=diamond,mark size=1.8pt]		table[x=intdofs, y=eff_index3]{\two};
\addplot[red,mark=o,mark size=1.8pt]		table[x=intdofs, y=eff_index2_union]{\two};
%
\end{loglogaxis}
\end{tikzpicture}
\hspace{-5pt}
\begin{tikzpicture}
\pgfplotstableread{tp2_multimesh_s07.dat}{\three}
\begin{loglogaxis}
[
title={\footnotesize $s = 0.7$},
width = 4.8cm, height = 4.8cm,						
xlabel={degrees of freedom (dof)}, 					
xlabel style={font=\fontsize{8pt}{10pt}\selectfont},		
ymajorgrids=true, xmajorgrids=true, grid style=dashed,	
xmin=200, xmax=40000,						
ymin = 2*10^(-5),	 ymax = 3.0,							
]
\addplot[black,mark=triangle,mark size=2pt]		table[x=intdofs, y=comperror3]{\three};
\addplot[blue,mark=square,mark size=1.4pt]		table[x=intdofs, y=comperror2_union]{\three};
\addplot[magenta,mark=diamond,mark size=1.8pt]		table[x=intdofs, y=eff_index3]{\three};
\addplot[red,mark=o,mark size=1.8pt]		table[x=intdofs, y=eff_index2_union]{\three};
%
\end{loglogaxis}
\end{tikzpicture}
\centerline{\ref{legend_tp2_1}}
\caption{Test case I (multimesh setting; $s \in \{0.3,\, 0.5,\, 0.7\}$):
evolution of global error estimates $\eta^m$, $\widetilde\eta^m$ and
the corresponding effectivity indices $\Theta^m$, $\widetilde\Theta^m$.}
\label{fig:test:case:I:estimates:effectivity}
\end{figure}

In Figure~\ref{fig:test:case:I:estimates:effectivity}, for each $s \in \{0.3,\, 0.5,\, 0.7\}$,
we visualize the evolution of global error estimates $\eta^m$, $\widetilde\eta^m$ and
the corresponding effectivity indices $\Theta^m$, $\widetilde\Theta^m$.
Table~\ref{tab:test:case:I:rates} gives a more quantitative representation of these results;
here, we report the number of degrees of freedom in the final approximation $\widetilde u_\kappa^{m_*} \in \widetilde{\mathcal S}^{m_*,1}$
generated by Algorithm~\ref{alg:algorithm_outline}
and the associated error estimates $\eta^{m_*}$, $\widetilde\eta^{m_*}$,
the number of degrees of freedom in the reference solution $u_{\text{ref}} \in {\mathcal S}^{\text{ref}}$
and the corresponding global error estimate $\widetilde\eta^{\text{ref}}$ computed using~\eqref{eq:union_mesh_fractional_estimate_def},
as well as the decay rates (with respect to the total number of degrees of freedom)
for $\eta^{m}$ and $\widetilde\eta^{m}$.
The decay rates in this and other test cases
are calculated from a linear regression fit on the corresponding values of error estimates
for the last~15 iterations of the adaptive loop, in order to avoid any pre-asymptotic regime to affect the results.
%
We also include in Table~\ref{tab:test:case:I:rates} the theoretical convergence rates
for uniformly refined meshes as predicted by~\cite[Theorem 4.3]{2015_bonito}.

We note that the dimension of the reference finite element space ${\mathcal S}^{\text{ref}}$ is at least 32 times bigger than
the number of degrees of freedom in the final approximation $\widetilde u_\kappa^{m_*}$ generated by the adaptive algorithm
(see Table~\ref{tab:test:case:I:rates}).
\rev{What is more important though is that}
the error estimate $\widetilde\eta^{\text{ref}}$ for the reference solution $u_{\text{ref}}$
is an order of magnitude smaller than the error estimate \rev{$\widetilde\eta^{m_*}$}.
This justifies the use of $u_{\text{ref}}$ as a proxy for the true solution to problem~\eqref{eq:fractional_laplacian_strong}
when calculating the effectivity indices in~\eqref{eq:effectivity:indices}.
The plots in Figure~\ref{fig:test:case:I:estimates:effectivity} show that
the global error estimation strategy based on the union mesh provides a more effective and robust way to estimate the error in
multimesh approximations than the strategy based on the triangle inequality.
Indeed, the effectivity indices $\widetilde\Theta^m$ (for the union mesh-based strategy)
are only slightly less than unity (which is typical for hierarchical error estimates) and,
more importantly, vary insignificantly across iterations,
whereas the effectivity indices $\Theta^m$ (for the triangle inequality-based strategy) exhibit a mild growth as iterations progress.
This results in deterioration of the convergence rate for the error estimates $\eta^m$
compared to that of the error estimates $\widetilde\eta^m$.
The convergence rates for the error estimates $\widetilde\eta^m$
are suboptimal (i.e., less than~1) and they improve as the fractional power $s$ gets closer to~1.
These rates are higher than those predicted by~\cite[Theorem 4.3]{2015_bonito}
for uniformly refined meshes, particularly for smaller $s$.

\begin{table}[!th]
  \caption{Test case I (multimesh setting; $s \in \{0.3,\, 0.5,\, 0.7\}$):
  the number of degrees of freedom and the error estimates for the final approximation
  $\widetilde u_\kappa^{m_*} \,{\in}\, \widetilde{\mathcal S}^{m_*,1}$
  generated by Algorithm~\ref{alg:algorithm_outline} and for the reference solution $u_{\text{ref}} \,{\in}\, {\mathcal S}^{\text{ref}}$,
  as well as the decay rates for two global error estimates $\eta^m$,~$\widetilde\eta^m$
  and the theoret\-ical convergence rates for uniform mesh
  refinement 
  (see~\cite[Theorem 4.3]{2015_bonito}).
  }
  \label{tab:test:case:I:rates}
  \begin{center}
    \begin{tabular}{lccc}
\toprule
$s$ & 0.3 & 0.5 & 0.7 \\
\midrule
$\text{dim}\big( \widetilde{\mathcal S}^{m_*,1} \big)$ & 9,485 & 13,365 & 34,781 \\
$\eta^{m_*}$ & 3.4e-03 & 8.2e-04 & 1.4e-04 \\
$\widetilde\eta^{m_*}$ & 1.8e-03 & 3.8e-04 & 6.1e-05 \\
\midrule
$\text{dim}\big( {\mathcal S}^{\text{ref}} \big)$ & 1,052,507 & 1,046,529 & 1,132,624 \\
$\widetilde\eta^{\text{ref}}$ & 1.7e-04 & 3.1e-05 & 6.0e-06 \\
\midrule
decay rate for $\eta^m$ & 0.70 & 0.81 & 0.90 \\
decay rate for $\widetilde \eta^m$ & 0.82 & 0.90 & 0.96 \\
decay rate for uniform refinement & 0.55 & 0.75 & 0.95 \\
\bottomrule
\end{tabular}

  \end{center}
\end{table}

\subsection{Test case II: discontinuous right-hand side function}\label{sec:test:case:II}

In this test case, we set $\Omega = (0,1)^2$ and
solve the model problem~\eqref{eq:fractional_laplacian_strong} with
discontinuous right-hand side function defined as
\begin{equation*}
  f(x,y) =
  \begin{cases}
    -1 & \text{if } (x,y) \in Q_1 \cup Q_2,\\
    1 & \text{otherwise},
  \end{cases}
\end{equation*}
where $Q_1$ and $Q_2$ are the subsets of $\Omega$ bounded by two quarter-circles of radius 0.6 centered at
$(0,0)$ and $(1,1)$, respectively.
The discontinuity of $f$ leads to sharp gradients in the solution $u$ along the interfaces,
in addition to $u$ exhibiting boundary layers.
Furthermore, the presence of curved interfaces in the definition of $f$
makes the use of a priori constructed graded meshes less straightforward than, e.g., in test case~I;
this further motivates the application of an adaptive mesh refinement algorithm.

In our first experiment for this test case, for $s \in \{0.1,\, 0.3,\, 0.5,\, 0.7,\, 0.9\}$,
we run the adaptive multimesh refinement algorithm
(Algorithm~\ref{alg:algorithm_outline})
as well as its single-mesh version developed in~\cite{2023_bulle}.
In each run, we start with the initial uniform mesh $\T^0$ of 512 right-angled triangles
(see Figure~\ref{fig:test:case:II:meshes}).
For each $s$, the number $N = N(\kappa)$ of parametric problems as well as the stopping tolerance $\tol$
are shown in Table~\ref{tab:test:case:II:rates}.

\begin{figure}[t!]
\begin{tikzpicture}
\pgfplotstableread{tp4_multimesh_s03.dat}{\one}
\pgfplotstableread{tp4_single-mesh_s03.dat}{\two}
\begin{loglogaxis}
[
title={\footnotesize $s = 0.3$},
width = 4.8cm, height = 4.8cm,						
xlabel={degrees of freedom (dof)}, 					
xlabel style={font=\fontsize{8pt}{10pt}\selectfont},		
ylabel={estimated error},				
ylabel style={font=\fontsize{8pt}{10pt}\selectfont}, 		
ymajorgrids=true, xmajorgrids=true, grid style=dashed,	
xmin=200, xmax=230000,						
ymin = 7*10^(-5),	 ymax = 2.7*10^(-2),							
legend columns=3,
legend entries={
$\widetilde\eta^m$ (single-mesh)\qquad\qquad,
$\eta^m$ (multimesh)\qquad\qquad,
$\widetilde\eta^m$ (multimesh)\qquad\qquad
},
legend to name=legend_tp4_1,
legend style={legend cell align=left, row sep=5pt, column sep=2pt, fill=none, draw=none, font={\fontsize{9pt}{12pt}\selectfont}},
]
\addplot[red,mark=o,mark size=1.8pt]		table[x=intdofs, y=comperror2]{\two};
\addplot[black,mark=triangle,mark size=2pt]		table[x=intdofs, y=comperror3]{\one};
\addplot[blue,mark=square,mark size=1.4pt]		table[x=intdofs, y=comperror2_union]{\one};
%
\end{loglogaxis}
\end{tikzpicture}
\hspace{-5pt}
\begin{tikzpicture}
\pgfplotstableread{tp4_multimesh_s05.dat}{\one}
\pgfplotstableread{tp4_single-mesh_s05.dat}{\two}
\begin{loglogaxis}
[
title={\footnotesize $s = 0.5$},
width = 4.8cm, height = 4.8cm,						
xlabel={degrees of freedom (dof)}, 					
xlabel style={font=\fontsize{8pt}{10pt}\selectfont},		
ymajorgrids=true, xmajorgrids=true, grid style=dashed,	
xmin=200, xmax=155000,						
ymin = 1.5*10^(-5),	 ymax = 8*10^(-3),							
]
\addplot[red,mark=o,mark size=1.8pt]		table[x=intdofs, y=comperror2]{\two};
\addplot[black,mark=triangle,mark size=2pt]		table[x=intdofs, y=comperror3]{\one};
\addplot[blue,mark=square,mark size=1.4pt]		table[x=intdofs, y=comperror2_union]{\one};
%
\end{loglogaxis}
\end{tikzpicture}
\hspace{-5pt}
\begin{tikzpicture}
\pgfplotstableread{tp4_multimesh_s07.dat}{\one}
\pgfplotstableread{tp4_single-mesh_s07.dat}{\two}
\begin{loglogaxis}
[
title={\footnotesize $s = 0.7$},
width = 4.8cm, height = 4.8cm,						
xlabel={degrees of freedom (dof)}, 					
xlabel style={font=\fontsize{8pt}{10pt}\selectfont},		
ymajorgrids=true, xmajorgrids=true, grid style=dashed,	
xmin=200, xmax=150000,						
ymin = 3*10^(-6),	 ymax = 2.3*10^(-3),							
]
\addplot[red,mark=o,mark size=1.8pt]		table[x=intdofs, y=comperror2]{\two};
\addplot[black,mark=triangle,mark size=2pt]		table[x=intdofs, y=comperror3]{\one};
\addplot[blue,mark=square,mark size=1.4pt]		table[x=intdofs, y=comperror2_union]{\one};
%
\end{loglogaxis}
\end{tikzpicture}
\centerline{\ref{legend_tp4_1}}
\caption{Test case II: evolution of global error estimates in the single-mesh setting 
($\widetilde\eta^m$ defined by~\eqref{eq:union_mesh_fractional_estimator1}) and
in the multimesh setting ($\eta^m$ and $\widetilde\eta^m$ defined
by~\eqref{eq:triangular_inequality_fractional_estimator}
and~\eqref{eq:union_mesh_fractional_estimate_def}, respectively) for $s \in \{0.3,\, 0.5,\, 0.7\}$.}
\label{fig:test:case:II:multimesh:vs:singlemesh}
\end{figure}
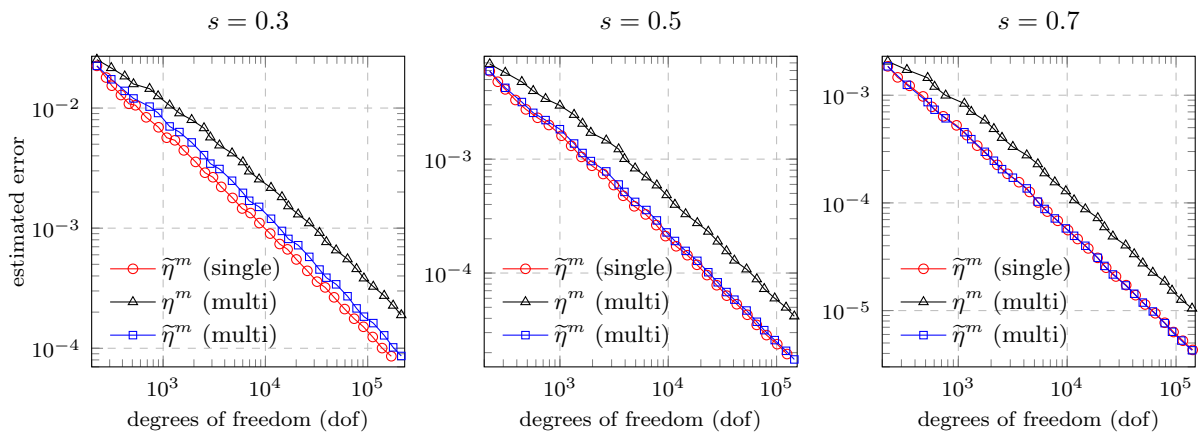

%
In Figure~\ref{fig:test:case:II:multimesh:vs:singlemesh},
for $s \in \{0.3,\, 0.5,\, 0.7\}$,
we report the evolution of three global a posteriori error estimates:
$\widetilde\eta^m$ defined by~\eqref{eq:union_mesh_fractional_estimator1} in the single-mesh setting
as well as
$\eta^m$ and $\widetilde\eta^m$ defined, respectively, by~\eqref{eq:triangular_inequality_fractional_estimator}
and~\eqref{eq:union_mesh_fractional_estimate_def} for multimesh discretizations.
%
%
The decay rates for these three global error estimates
%
as well as the theoretical convergence rates predicted in~\cite[Theorem 4.3]{2015_bonito}
for uniform mesh refinement are shown in Table~\ref{tab:test:case:II:rates}
for $s \in \{0.1,\, 0.3,\, 0.5,\, 0.7,\, 0.9\}$.

We make the following conclusions by looking at Figure~\ref{fig:test:case:II:multimesh:vs:singlemesh}
and Table~\ref{tab:test:case:II:rates}.
Firstly, in the multimesh setting, similarly to test case~I, we observe a slower decay of error estimates $\eta^m$
compared to that of $\widetilde \eta^m$.
This is again due to the triangle inequality affecting the quality of the error estimation
in~\eqref{eq:true_error_bound_triangular_inequality_estimator}.
Secondly, while convergence rates for adaptive algorithms (in both the single-mesh and multimesh settings)
are suboptimal 
for $s \in \{0.1,\, 0.3,\, 0.5,\, 0.7\}$ as expected, these rates exceed
the theoretically predicted and experimentally observed
rates for uniformly refined meshes (in the single-mesh setting),
particularly for smaller values of $s$ (cf.~\cite{2015_bonito}).
This demonstrates the advantage of adaptive mesh refinement algorithms in mitigating the singular behavior
of the solution $u$ along the discontinuity interfaces in $f$.
Again, this advantage is more pronounced for smaller values of $s$.
Finally, we only see a marginal improvement of convergence rates in the multimesh setting
compared to adaptive single-mesh discretizations, and this improvement diminishes as $s$ increases.
For $s=0.9$, both adaptive single-mesh and multimesh discretizations converge with essentially an optimal rate
$\mathcal{O}\big(\mathrm{dof}^{-1}\big)$---the predicted convergence rate
for uniformly refined (single-mesh) approximations for this value of $s$.

\begin{table}[!th]
  \caption{Test case II:
  the number $N$ of parametric problems,
  the stopping tolerance $\tol$, 
  the number $m_*$ of iterations as well as
  the decay rates for global error estimates in the single-mesh and multimesh settings,
  and the theoret\-ical convergence rates for uniform 
  refinement in the single-mesh setting (see~\cite[Theorem 4.3]{2015_bonito}).
  }
  \label{tab:test:case:II:rates}
  \begin{center}
    \begin{tabular}{lccccc}
\toprule
$s$ & 0.1 & 0.3 & 0.5 & 0.7 & 0.9 \\
\midrule
$N$ & 408 & 176 & 149 & 176 & 408 \\
$\tol$ & 1e-03 & 1e-04 & 2e-05 & 5e-06 & 1e-06 \\
\midrule
$m_*$ (adaptive single-mesh) & 33 & 32 & 29 & 30 & 32 \\
$m_*$ (adaptive multimesh) & 35 & 30 & 28 & 28 & 30 \\
\midrule
decay rate for $\widetilde \eta^m$ (adaptive single-mesh) & 0.66 & 0.86 & 0.93 & 0.97 & 0.99 \\
decay rate for $\eta^m$ (adaptive multimesh) & 0.64 & 0.82 & 0.90 & 0.93 & 0.96 \\
decay rate for $\widetilde \eta^m$ (adaptive multimesh) & 0.67 & 0.89 & 0.96 & 0.98 & 0.99 \\
decay rate for uniform refinement 
& 0.35 & 0.55 & 0.75 & 0.95 & 1.00 \\
\bottomrule
\end{tabular}

  \end{center}
\end{table}

\begin{figure}[!th]

  \begin{center}
      \includegraphics[height=3.4cm]{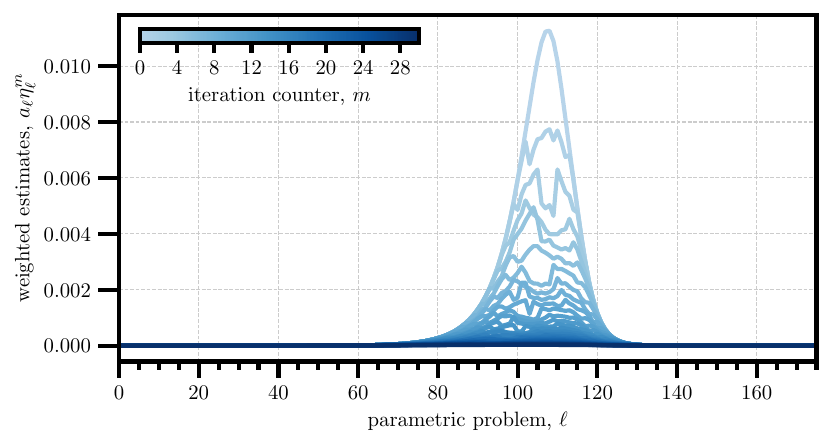}
      \includegraphics[height=3.4cm]{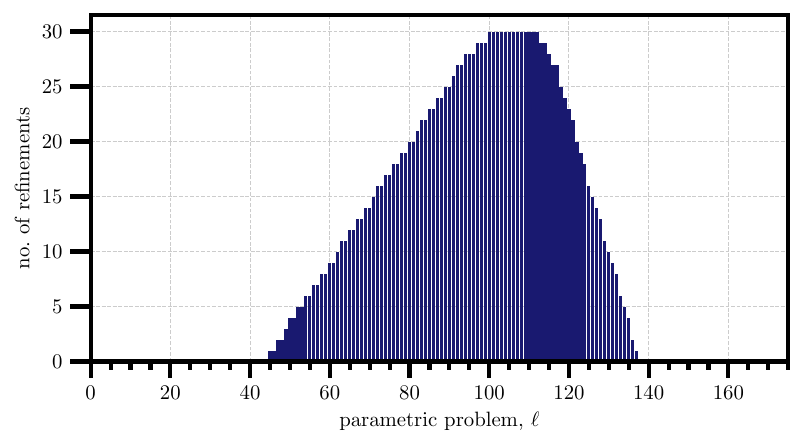}\\
      \includegraphics[height=3.4cm]{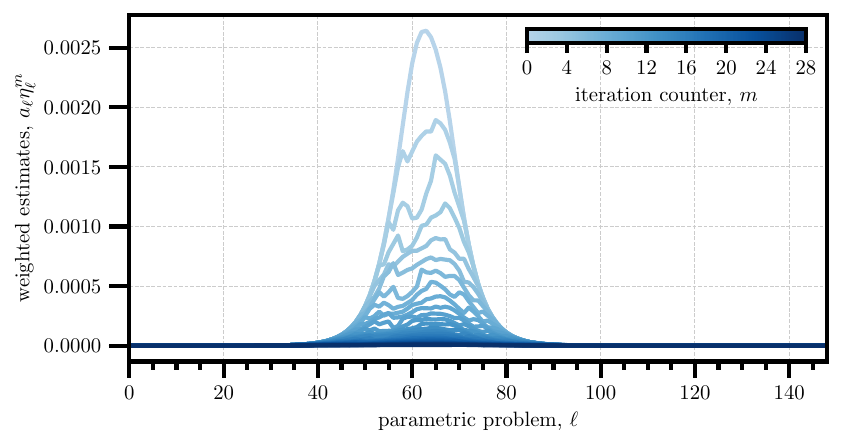}
      \includegraphics[height=3.4cm]{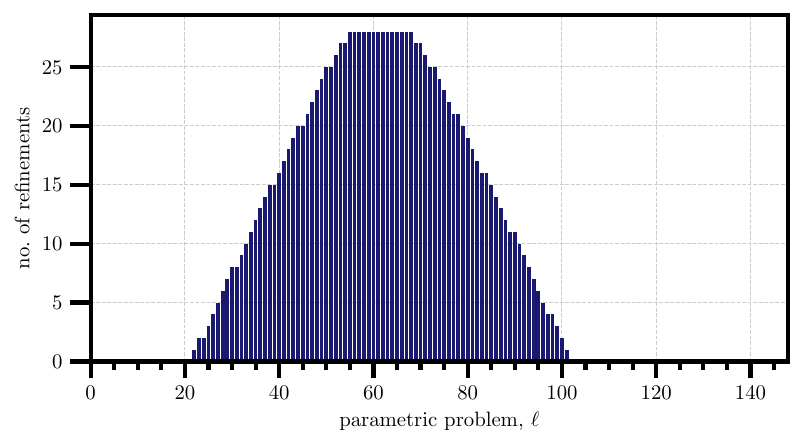}\\
      \includegraphics[height=3.4cm]{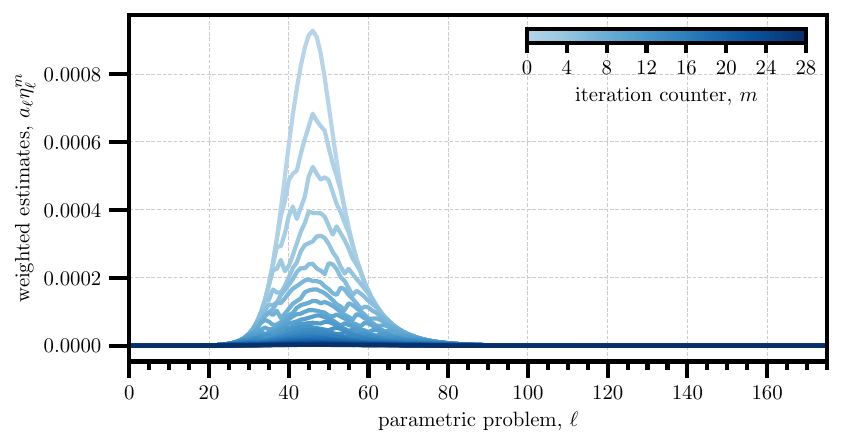}
      \includegraphics[height=3.4cm]{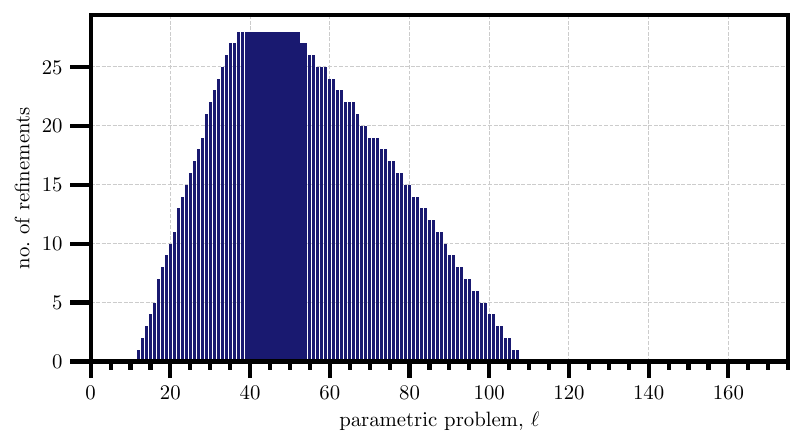}
  \end{center}
  \caption{Test case II (multimesh setting;
  $s=0.3$ (top row), $s=0.5$ (middle row), and $s=0.7$ (bottom row)):
  the number of refinements for each mesh $\T_l$ ($l = 1,\ldots,N$) (right column) and
  the corresponding weighted error estimates $a_l \eta_l^m$
  (see~\eqref{eq:global_parametric_estimator}) (left column).
  In the left column, different shades represent different iterations~$m$ of the adaptive loop.
  }
  \label{fig:test:case:II:param:refinement:estimators}
\end{figure}
\begin{figure}[!th]
  \begin{center}
      \includegraphics*[height=4.5cm, trim = 52 0 72 0, clip]{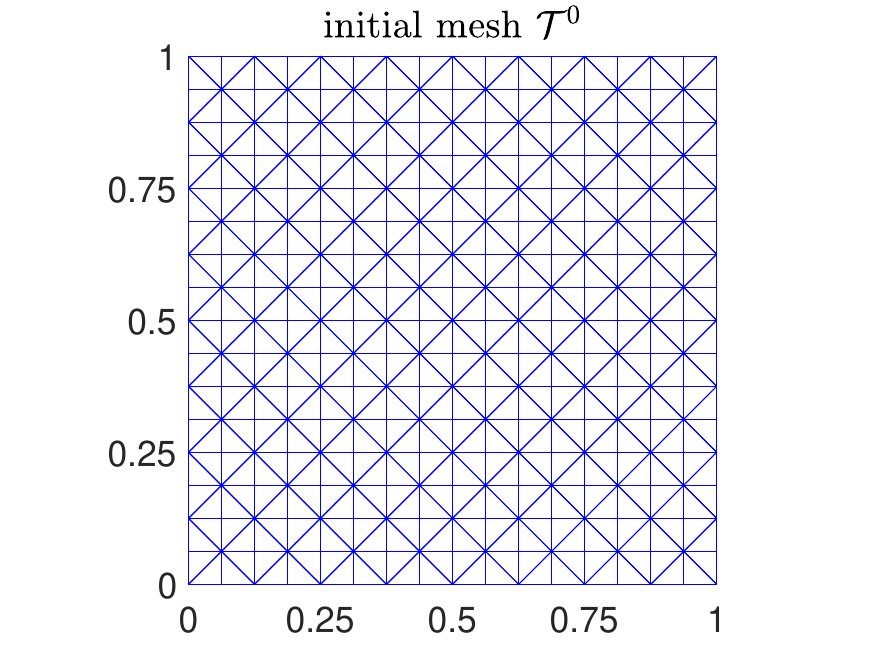}
      \
      \includegraphics*[height=4.5cm, trim = 85 25 72 0, clip]{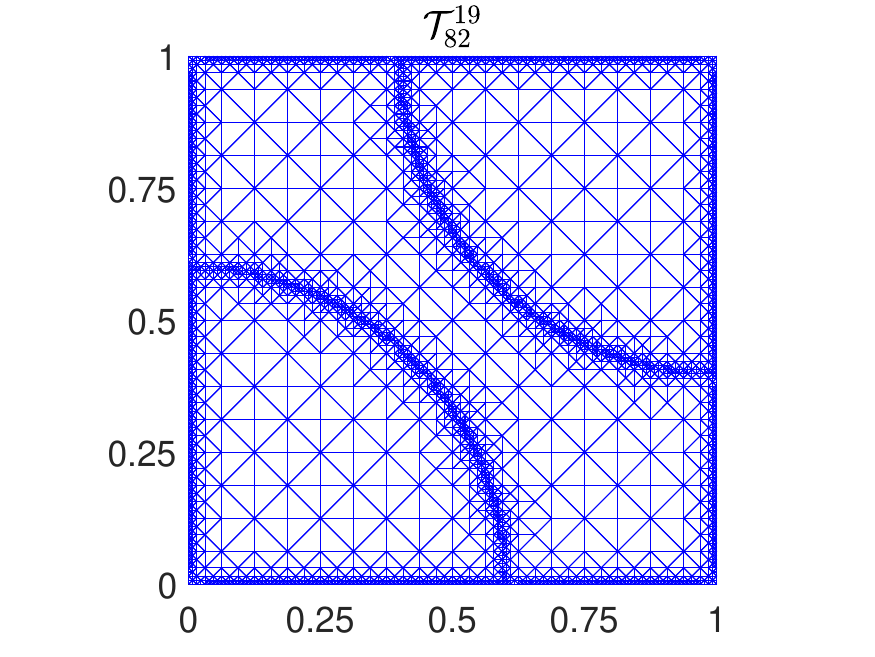}
      \
      \includegraphics*[height=4.5cm, trim = 85 25 72 0, clip]{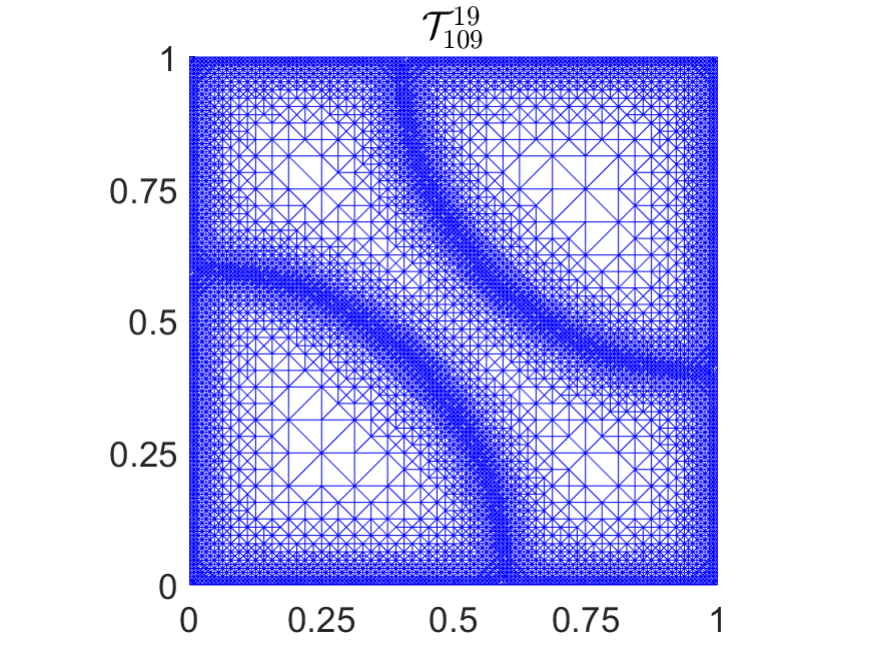}
      \\[7pt]
      \includegraphics*[height=4.5cm, trim = 85 25 72 0, clip]{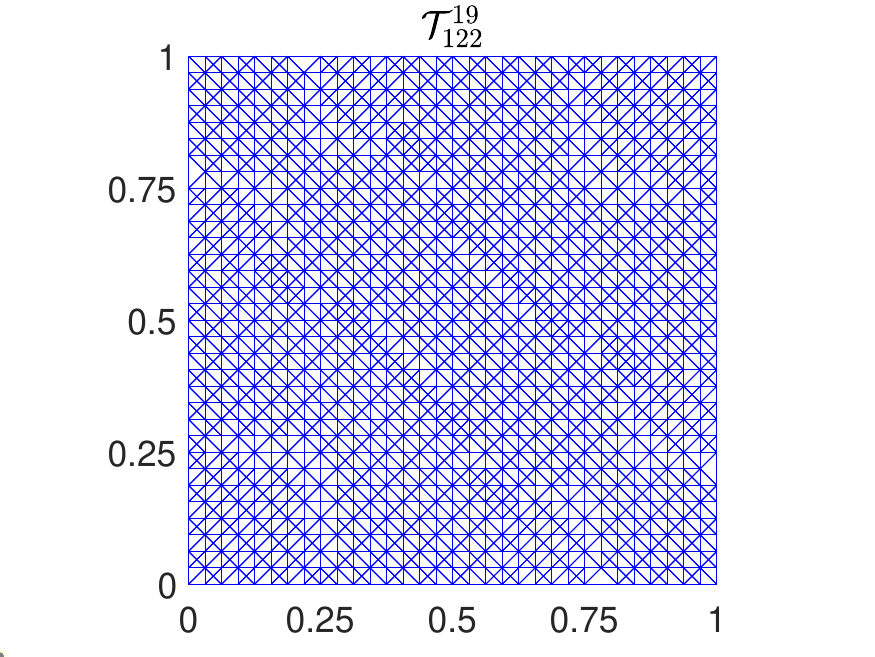}
      \
      \includegraphics*[height=4.5cm, trim = 85 25 72 0, clip]{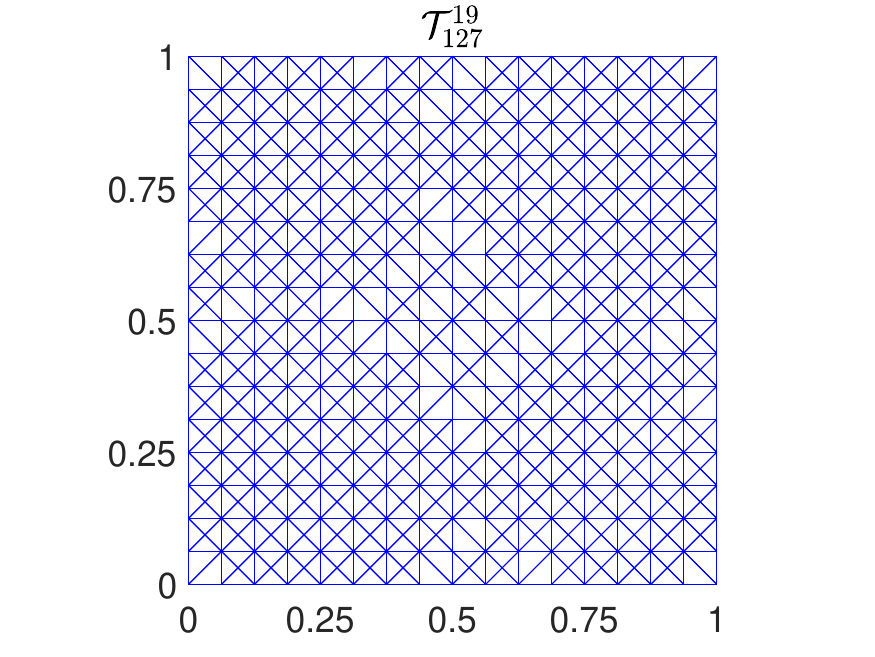}
      \
      \includegraphics*[height=4.5cm, trim = 86 25 72 0, clip]{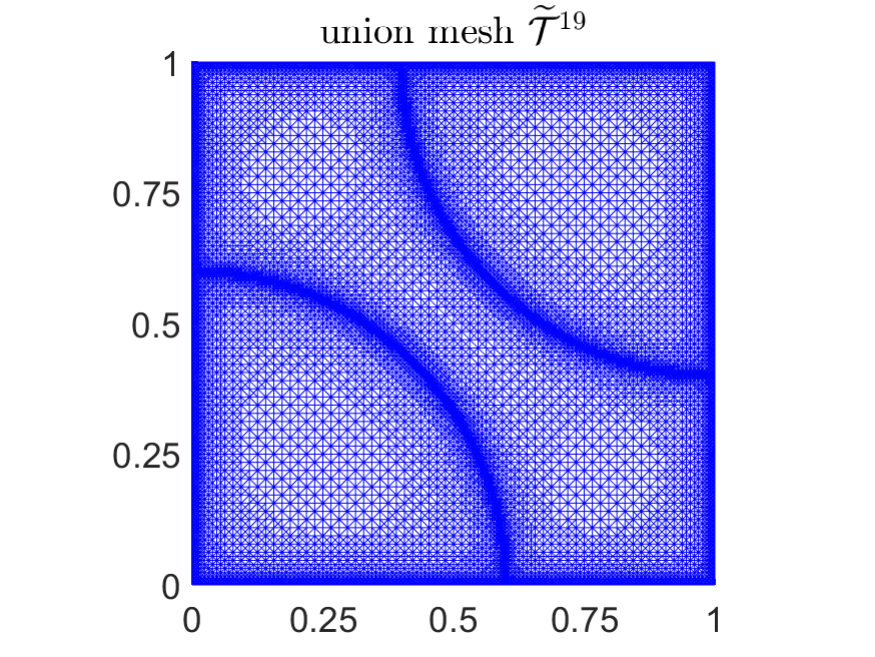}
    \end{center}
  \caption{Test case II (multimesh setting; $s=0.3$):
  the initial coarse mesh $\T^0$;
  finite element meshes $\T_l^m$ generated by Algorithm~\ref{alg:algorithm_outline} ($m=19$; $l=82,\, 109,\, 122,\, 127$);
  the union mesh $\widetilde\T^{19}$.
  }
  \label{fig:test:case:II:meshes}
\end{figure}

In Figure~\ref{fig:test:case:II:param:refinement:estimators}, for
$s \in \{0.3,\, 0.5,\, 0.7\}$,
we plot the values of weighted (global) estimates $a_l \eta_l^m$
(see~\eqref{eq:global_parametric_estimator})
for each parametric problem (i.e., for $l=1\ldots N$) and
across all iterations of the adaptive loop (i.e., for $m=1,\ldots,m_*$);
we also plot the number of times the mesh $\T_l$ for the $l$-th parametric problem is refined.
These plots show that our multimesh refinement algorithm is successful in detecting the parametric problems where the weighted error estimates are large,
and it concentrates most of refinements on the corresponding meshes.
Furthermore, we observe that the underlying meshes are not refined for
about half of parametric problems
(specifically, $47\%$ for $s=0.3$, $46\%$ for $s=0.5$, and $45\%$ for $s=0.7$);
therefore, these parametric problems are only solved once on the coarsest mesh.
In addition, we can see from the 
plots in the left column of Figure~\ref{fig:test:case:II:param:refinement:estimators}
that as iterations progress, the weighted error estimates become more and more even in magnitude 
across parametric problems.
This is ensured by D{\"o}rfler marking~\eqref{eq:marking}.


Figure~\ref{fig:test:case:II:meshes} depicts the initial coarse mesh $\T^0$,
the finite element meshes generated at the 19th iteration the adaptive multimesh refinement algorithm
for some parametric problems, and the corresponding union mesh for test case~II with $s=0.3$.

Let us now compare the computational costs of running the adaptive algorithm
in the single-mesh and multimesh settings.
First, we assess the costs in terms of the number of degrees of freedom in computed approximations.
For $l \in \{1,\ldots,N\}$, we define $\mathrm{cost}_l^m$ as the number of degrees of freedom
in the finite element formulation~\eqref{eq:discrete_parametric_problems} of the corresponding parametric problem
(in other words, $\mathrm{cost}_l^m = \mathrm{dim} \big( \mathcal S^{m,1}_l \big)$).
Then, the overall computational cost at iteration $m \in \N_0$ of the adaptive algorithm
is the sum of $\mathrm{cost}_l^m$ across all the parametric problems that are actually solved at this iteration
(for $m \in \N$, this excludes those parametric problems for which the underlying
meshes have not been refined at iteration $m-1$).
Thus, recalling that the meshes for all parametric problems are initialized by the coarse mesh $\T^0$,
the total computational cost at iteration $m \in \N_0$ is given~by
\begin{equation}
  \mlvar{\mathrm{totcost}}^m :=
  \begin{cases}
    \displaystyle 
                         N \times \mathrm{cost}_1^0    & \text{if } m=0,\\
    \displaystyle \sum_{l\in \mathcal L} \mathrm{cost}_l^m & \text{otherwise},
  \end{cases}
  \label{eq:totcost_multi_mesh}
\end{equation}
where $\mathcal L \subseteq \{1,\ldots,N\}$ determines the subset of parametric problems
for which the underlying meshes have been refined at iteration~$m-1$.
Note that in the single-mesh setting 
formula~\eqref{eq:totcost_multi_mesh} simplifies to
\begin{equation*}
  \mlvar{\mathrm{totcost}}^m = N \times \mathrm{cost}_1^m,\quad
  m \in \N_0.
\end{equation*}
The cumulative cost at the $m$-th iteration of the adaptive algorithm is defined as follows:
\begin{equation}
  \mathrm{cumcost}^m := \sum_{j=1}^m \mlvar{\mathrm{totcost}}^j.
  \label{eq:cumulative_cost}
\end{equation}
Figure~\ref{fig:test:case:II:cost:comparison} 
shows the growth of cumulative costs for computing adaptively refined approximations.
Here, we plot the cumulative costs (as defined in~\eqref{eq:cumulative_cost})
against global a posteriori error estimates $\widetilde\eta^m$
in the single-mesh and multimesh settings
(see~\eqref{eq:union_mesh_fractional_estimator1} and~\eqref{eq:union_mesh_fractional_estimate_def}, respectively).
The plots show that computing an adaptive multimesh approximation can be up to 10 times cheaper
(in terms of the cumulative number of degrees of freedom)
than computing the adaptive single-mesh approximation to the same error tolerance.

\begin{figure}[!th]
\begin{tikzpicture}
\pgfplotstableread{tp4_multimesh_s01.dat}{\one}
\pgfplotstableread{tp4_single-mesh_s01.dat}{\two}
\begin{loglogaxis}
[
title={\footnotesize $s = 0.1$},
width = 4.8cm, height = 4.8cm,						
xlabel={estimated error}, 					
xlabel style={font=\fontsize{8pt}{10pt}\selectfont},		
ylabel={cumulative cost},				
ylabel style={font=\fontsize{8pt}{10pt}\selectfont}, 		
ymajorgrids=true, xmajorgrids=true, grid style=dashed,	
xmin = 7*10^(-4),	 xmax = 1.2*10^(-1),				
ymin=7*10^4, ymax=1.3*10^9,						
legend style={legend pos=south west, legend cell align=left, fill=none}
]
\addplot[red,mark=o,mark size=1.8pt]		table[x=comperror2, y=cumul_cost]{\two};
\addplot[blue,mark=square,mark size=1.4pt]		table[x=comperror2_union, y=cumul_cost]{\one};
%
\end{loglogaxis}
\end{tikzpicture}
\hspace{-5pt}
\begin{tikzpicture}
\pgfplotstableread{tp4_multimesh_s03.dat}{\one}
\pgfplotstableread{tp4_single-mesh_s03.dat}{\two}
\begin{loglogaxis}
[
title={\footnotesize $s = 0.3$},
width = 4.8cm, height = 4.8cm,						
xlabel={estimated error}, 					
xlabel style={font=\fontsize{8pt}{10pt}\selectfont},		
ymajorgrids=true, xmajorgrids=true, grid style=dashed,	
xmin = 7*10^(-5),	 xmax = 2.7*10^(-2),				
ymin=3.5*10^4, ymax=2.4*10^8,						
legend style={legend pos=south west, legend cell align=left, fill=none}
]
\addplot[red,mark=o,mark size=1.8pt]		table[x=comperror2, y=cumul_cost]{\two};
\addplot[blue,mark=square,mark size=1.4pt]		table[x=comperror2_union, y=cumul_cost]{\one};
%
\end{loglogaxis}
\end{tikzpicture}
\hspace{-5pt}
\begin{tikzpicture}
\pgfplotstableread{tp4_multimesh_s05.dat}{\one}
\pgfplotstableread{tp4_single-mesh_s05.dat}{\two}
\begin{loglogaxis}
[
title={\footnotesize $s = 0.5$},
width = 4.8cm, height = 4.8cm,						
xlabel={estimated error}, 					
xlabel style={font=\fontsize{8pt}{10pt}\selectfont},		
ymajorgrids=true, xmajorgrids=true, grid style=dashed,	
xmin = 1.4*10^(-5),	 xmax = 7*10^(-3),				
ymin=3.5*10^4, ymax=1.5*10^8,						
legend style={legend pos=south west, legend cell align=left, fill=none}
]
\addplot[red,mark=o,mark size=1.8pt]		table[x=comperror2, y=cumul_cost]{\two};
\addplot[blue,mark=square,mark size=1.4pt]		table[x=comperror2_union, y=cumul_cost]{\one};
%
\end{loglogaxis}
\end{tikzpicture}
\\
\begin{center}
\begin{tikzpicture}
\pgfplotstableread{tp4_multimesh_s07.dat}{\one}
\pgfplotstableread{tp4_single-mesh_s07.dat}{\two}
\begin{loglogaxis}
[
title={\footnotesize $s = 0.7$},
width = 4.8cm, height = 4.8cm,						
xlabel={estimated error}, 					
xlabel style={font=\fontsize{8pt}{10pt}\selectfont},		
ylabel={cumulative cost},				
ylabel style={font=\fontsize{8pt}{10pt}\selectfont}, 		
ymajorgrids=true, xmajorgrids=true, grid style=dashed,	
xmin = 3*10^(-6),	 xmax = 2.3*10^(-3),				
ymin=3.5*10^4, ymax=2*10^8,						
legend style={legend pos=south west, legend cell align=left, fill=none}
]
\addplot[red,mark=o,mark size=1.8pt]		table[x=comperror2, y=cumul_cost]{\two};
\addplot[blue,mark=square,mark size=1.4pt]		table[x=comperror2_union, y=cumul_cost]{\one};
%
\end{loglogaxis}
\end{tikzpicture}
\hspace{10pt}
\begin{tikzpicture}
\pgfplotstableread{tp4_multimesh_s09.dat}{\one}
\pgfplotstableread{tp4_single-mesh_s09.dat}{\two}
\begin{loglogaxis}
[
title={\footnotesize $s = 0.9$},
width = 4.8cm, height = 4.8cm,						
xlabel={estimated error}, 					
xlabel style={font=\fontsize{8pt}{10pt}\selectfont},		
ymajorgrids=true, xmajorgrids=true, grid style=dashed,	
xmin = 6.4*10^(-7),	 xmax = 9*10^(-4),				
ymin=7*10^4, ymax=8*10^8,						
legend style={legend pos=south west, legend cell align=left, fill=none}
]
\addplot[red,mark=o,mark size=1.8pt]		table[x=comperror2, y=cumul_cost]{\two};
\addplot[blue,mark=square,mark size=1.4pt]		table[x=comperror2_union, y=cumul_cost]{\one};
%
\end{loglogaxis}
\end{tikzpicture}
\end{center}
\caption{Test case II: cumulative costs (as defined in~\eqref{eq:cumulative_cost})
plotted against the global error estimates $\widetilde\eta^m$ at each iteration of adaptive loop
in the single-mesh (circular markers) and multimesh (square markers) settings.
}
\label{fig:test:case:II:cost:comparison}
\end{figure}

While the comparison based on the cumulative number of degrees of freedom 
provides a good estimate of the computational effort involved in running each version of the adaptive algorithm,
it does not take into account some non-negligible costs associated with running the multimesh version;
these include, e.g., the costs of assembling the union mesh $\widetilde\T^{m}$,
computing the (fully discrete) solution $\widetilde u_\kappa^{m}$ and the associated error estimate $\widetilde\eta^{m}$.
For a fair comparison of the two versions in terms of actual computational times,
let us now perform the second experiment for test case~II.
For $s \in \{0.1,\, 0.3,\, 0.5,\, 0.7,\, 0.9\}$, we run both versions of the adaptive algorithm for~31 iterations
(i.e., rather than setting a stopping tolerance, we set $m_* = 30$ in Algorithm~\ref{alg:algorithm_outline}).
In the single-mesh version, we compute the global error estimate $\widetilde \eta^m$
given by~\eqref{eq:union_mesh_fractional_estimator1} at each iteration
(recall that the union mesh does not need to be assembled in this case).
In the multimesh version, we make a slight modification to Algorithm~\ref{alg:algorithm_outline} so that
the global error estimate $\eta^m$ is computed at each iteration
(this is a cheaper global error estimate based on the triangle inequality,
see~\eqref{eq:true_error_bound_triangular_inequality_estimator}),
whereas a more expensive global error estimate based on the union mesh (see~\eqref{eq:union_mesh_fractional_estimate_def})
is only computed at the final iteration (i.e., we set $k = 30$ in Algorithm~\ref{alg:algorithm_outline}
and therefore, the union mesh $\widetilde\T^{30}$ as well as
the associated fully discrete approximation $\widetilde u_\kappa^{30}$ and the error estimate $\widetilde\eta^{30}$
are only computed once).

In Table~\ref{tab:test:case:II:times}, for each version, we record the computed global error estimates at the final iteration
as well as the total CPU times of running the algorithm and, in the multimesh setting,
the CPU times of running individual blocks\footnote{The computations were performed using MATLAB version R2024a on
a workstation with an Intel Xeon W-3265 2.70GHz CPU and 128GB of RAM.}.
Firstly, we note that for each $s \in \{0.1,\, 0.3,\, 0.5,\, 0.7,\, 0.9\}$,
the error estimate $\widetilde\eta^{30}$ in the multimesh setting
is close to the error estimate $\widetilde\eta^{30}$ in the single-mesh setting
(in fact, this was the case across all iterations in our first experiment for this test case,
see the plots in Figure~\ref{fig:test:case:II:multimesh:vs:singlemesh}).
This justifies the premise of this experiment---running both versions for the same number of iterations.
Secondly, as expected, we observe that superior effectivity of the global error estimate $\widetilde\eta^{30}$
in the multimesh setting comes at significant computational costs:
the CPU time of computing $\widetilde\eta^{30}$ often doubles the CPU time of running all 31 iterations of the algorithm.
Nevertheless, even if the CPU times of computing $\widetilde\eta^{30}$ are included
into the total CPU times for the multimesh version,
the improvement against the CPU times for running the single-mesh version is by a factor of~2.7 on average.
If, on the other hand, only the cheaper global error estimates $\eta^{m}$ are computed in the multimesh setting,
then the multimesh version outperforms the single-mesh one (in terms of total CPU times) by a factor of~7.3 on average
(here, the total CPU times for the multimesh version include the costs of computing the union mesh~$\widetilde\T^{30}$
and the associated fully discrete approximation~$\widetilde u_\kappa^{30}$).

\begin{table}[!th]
  \caption{Test case II:
  final error estimates and CPU times (in seconds) of running the adaptive algorithm for 31 iterations
  in the single-mesh and multimesh settings.
  }
  \label{tab:test:case:II:times}
  \begin{center}
    {\small
\begin{tabular}{p{6.2cm}ccccc}
\toprule
$s$ & 0.1 & 0.3 & 0.5 & 0.7 & 0.9 \\
\midrule
\multicolumn{6}{c}
{\textbf{adaptive single-mesh algorithm running for 31 iterations ($m_* = 30$)}} \\
\midrule
error estimate $\widetilde\eta^{30}$
 & 1.4e-3 & 1.2e-4 & 1.6e-5 & 4.3e-6 & 1.4e-6 \\[3pt]
total CPU time for running the algorithm\hfill\break
(includes the computation of $\widetilde u_\kappa^{30}$ and $\widetilde \eta^{30}$)
 & 45061 & 12776 & 14099 & 13281 & 36277 \\
\midrule
\multicolumn{6}{c}
{\textbf{adaptive multimesh algorithm running for 31 iterations ($m_* = 30$)}} \\
\midrule
error estimate $\eta^{30}$
 & 2.7e-3 & 1.9e-4 & 2.8e-5 & 7.0e-6 & 2.4e-6 \\
error estimate $\widetilde\eta^{30}$
 & 1.8e-3 & 8.6e-5 & 1.1e-5 & 3.0e-6 & 9.7e-7 \\[3pt]
CPU time for running 31 iterations 
 & 3850 & 1385 & 1526 & 1860 & 3952 \\
CPU time for assembling 
$\widetilde\T^{30}$
 & 435 & 256 & 337 & 464 & 1080 \\
CPU time for computing 
$\widetilde u_\kappa^{30}$
 & 496 & 123 & 58 & 31 & 39 \\
CPU time for computing 
$\widetilde\eta^{30}$
 & 9864 & 3524 & 3026 & 3229 & 7317 \\
total CPU time 
(excl. computing $\widetilde\eta^{30}$)
 & 4781 & 1764 & 1921 & 2355 & 5071 \\
total CPU time 
(incl. computing $\widetilde\eta^{30}$)
 & 14645 & 5288 & 4947 & 5584 & 12388 \\[3pt]
Improvement factor (if $\widetilde\eta^{30}$ is not computed) &
 9.4 & 7.2 & 7.3 & 5.6 & 7.2 \\
Improvement factor (if $\widetilde\eta^{30}$ is computed) &
 3.1 & 2.4 & 2.9 & 2.4 & 2.9 \\
\bottomrule
\end{tabular}
}

  \end{center}
\end{table}

The performance gains achieved by the multimesh approach---both in terms of the overall number of degrees of freedom
and actual computational time---are of no surprise.
The marking strategy in Algorithm~\ref{alg:algorithm_outline} has been designed
to ensure that the algorithm does not over-refine the meshes associated with parametric problems
whose impact to the fully discrete rational approximation is small.
Furthermore, as discussed above (see Figure~\ref{fig:test:case:II:param:refinement:estimators}),
a large number of meshes for parametric problems are not actually refined at a given iteration
of Algorithm~\ref{alg:algorithm_outline}, meaning that the subroutines {\tt Solve} and {\tt Estimate}
are not run for these parametric problems at the next iteration.
Both these features of the multimesh approach contribute to significant reduction of computational costs.

\subsection{Test case III: L-shaped domain}\label{sec:test:case:III}

We now set  $f \equiv 1$ and look to solve the model problem~\eqref{eq:fractional_laplacian_strong}
on the L-shaped domain $\Omega = (-1,1)^2 \setminus (-1,0]^2$.
In this case, the solution $u$ to problem~\eqref{eq:fractional_laplacian_strong}
exhibits boundary layers 
as well as a geometric singularity at the domain's reentrant corner.
Therefore, the adaptive refinement strategy will need to capture
the interplay of these two different types of singular behavior in the solution~$u$.

\begin{figure}[!th]
    \begin{center}
        \includegraphics*[height=4.5cm, trim = 18 149 19 100, clip]{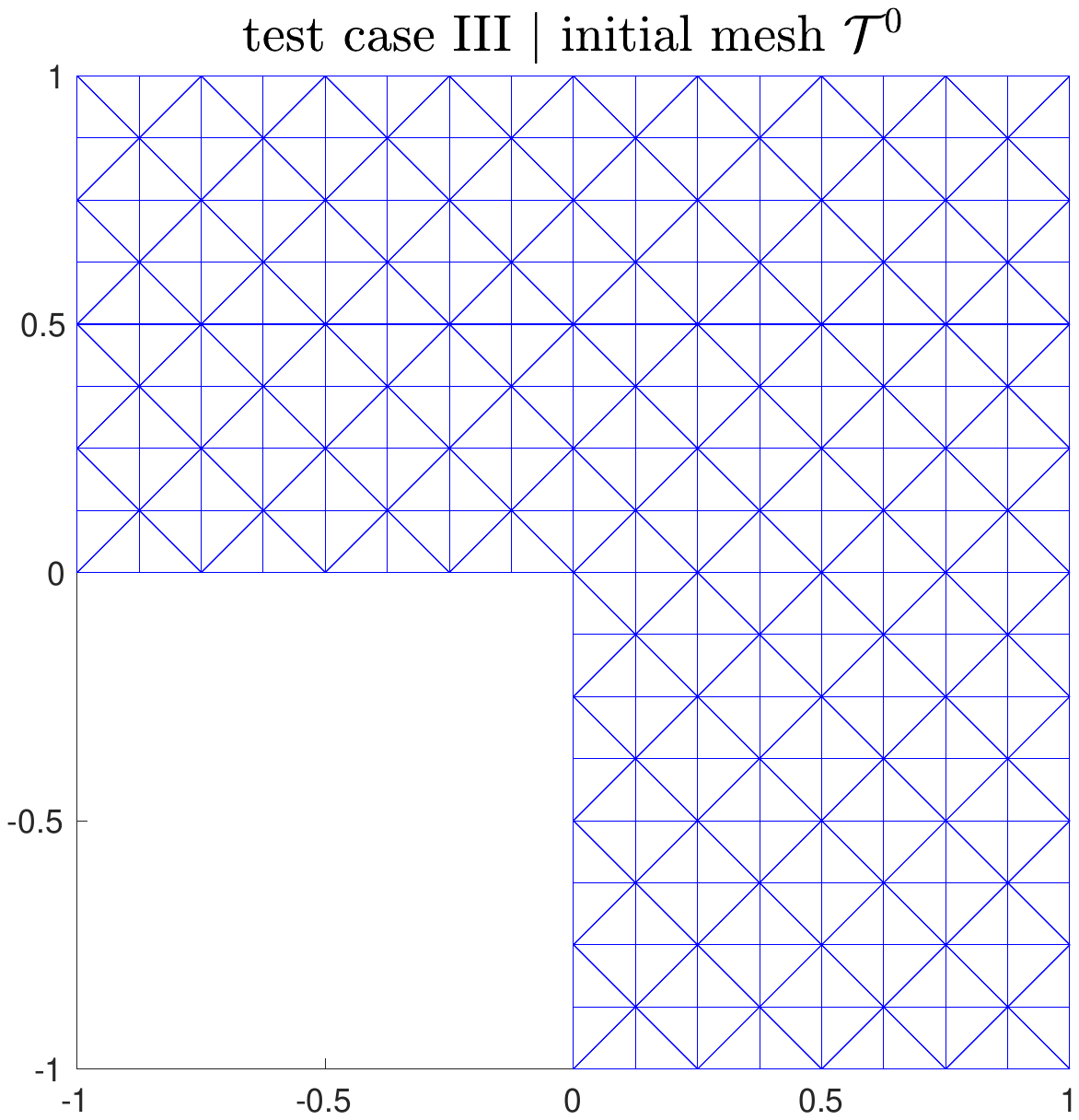}
        \
        \includegraphics*[height=4.5cm, trim = 50 165 19 100, clip]{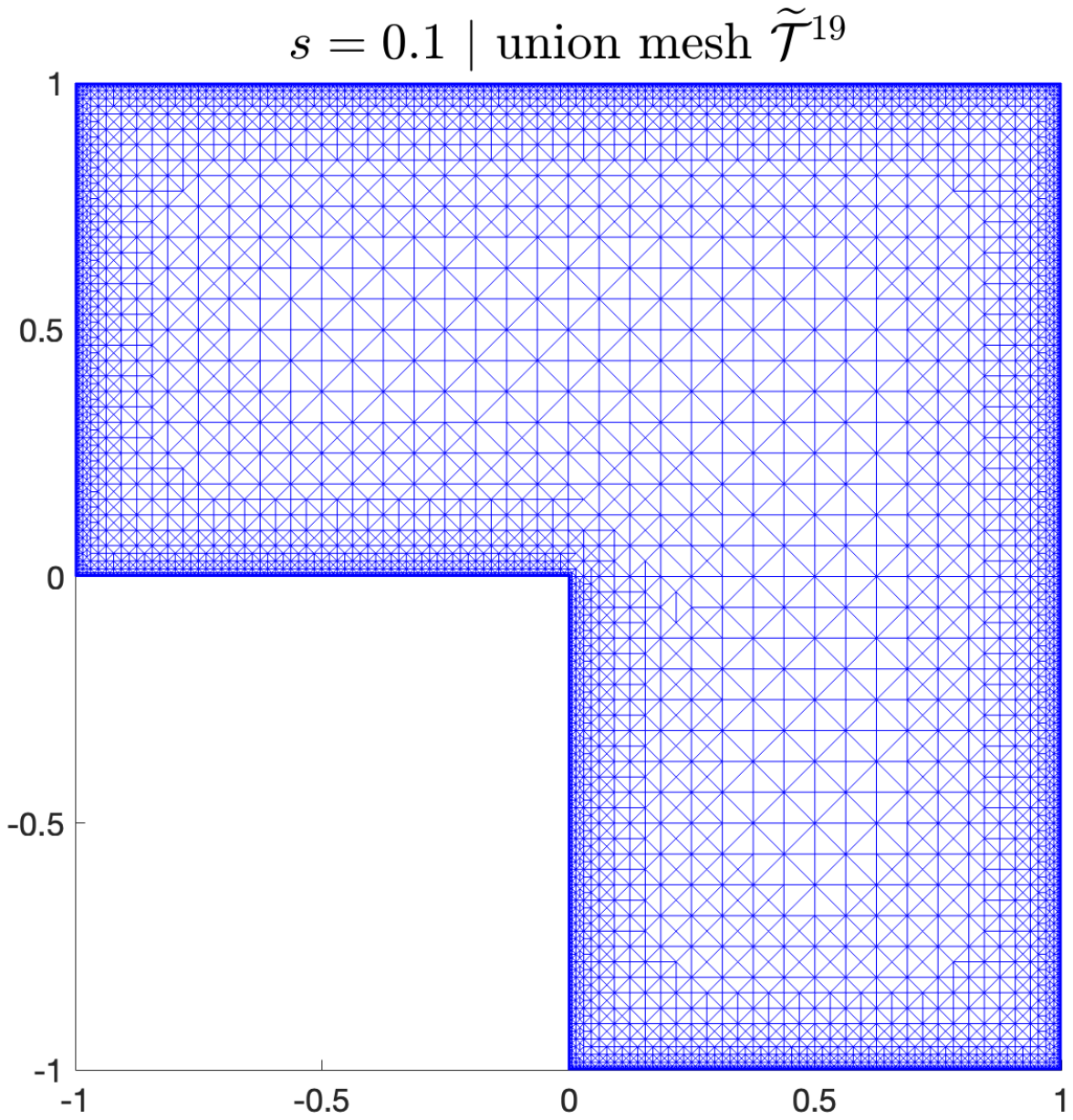}
        \
        \includegraphics*[height=4.5cm, trim = 50 165 19 100, clip]{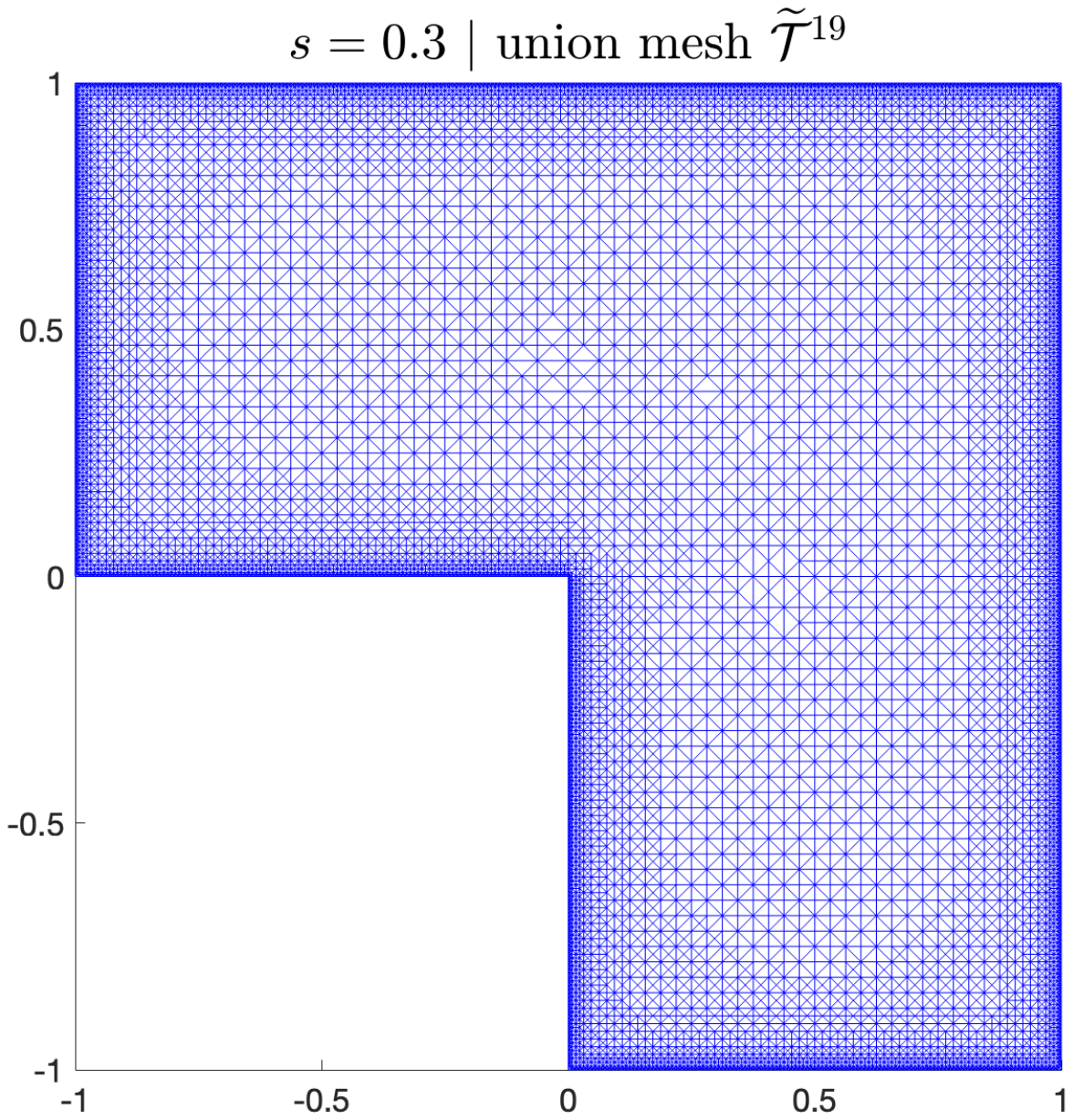}
        \\[7pt]
        \includegraphics*[height=4.5cm, trim = 50 165 19 100, clip]{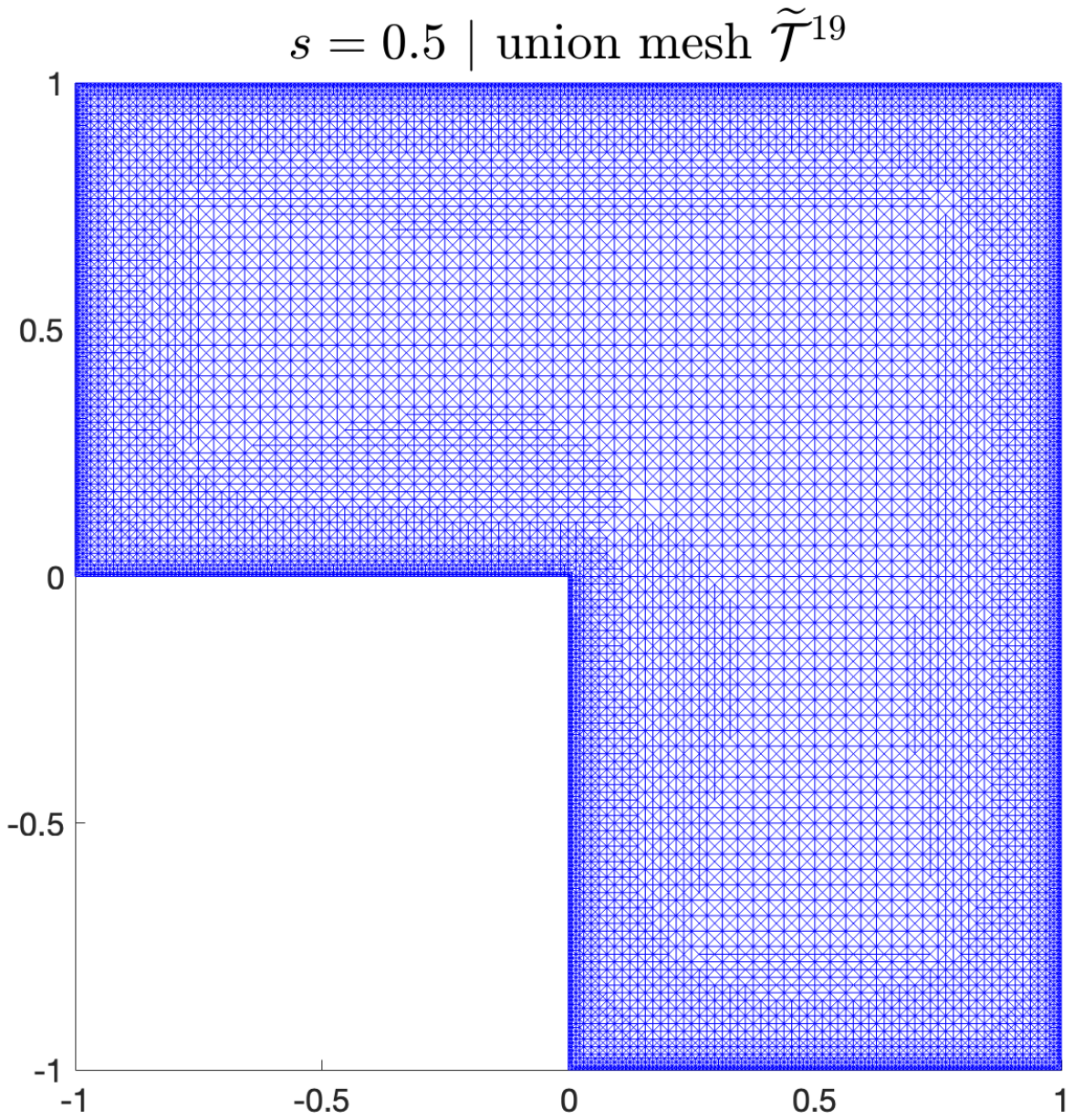}
        \
        \includegraphics*[height=4.5cm, trim = 50 165 19 100, clip]{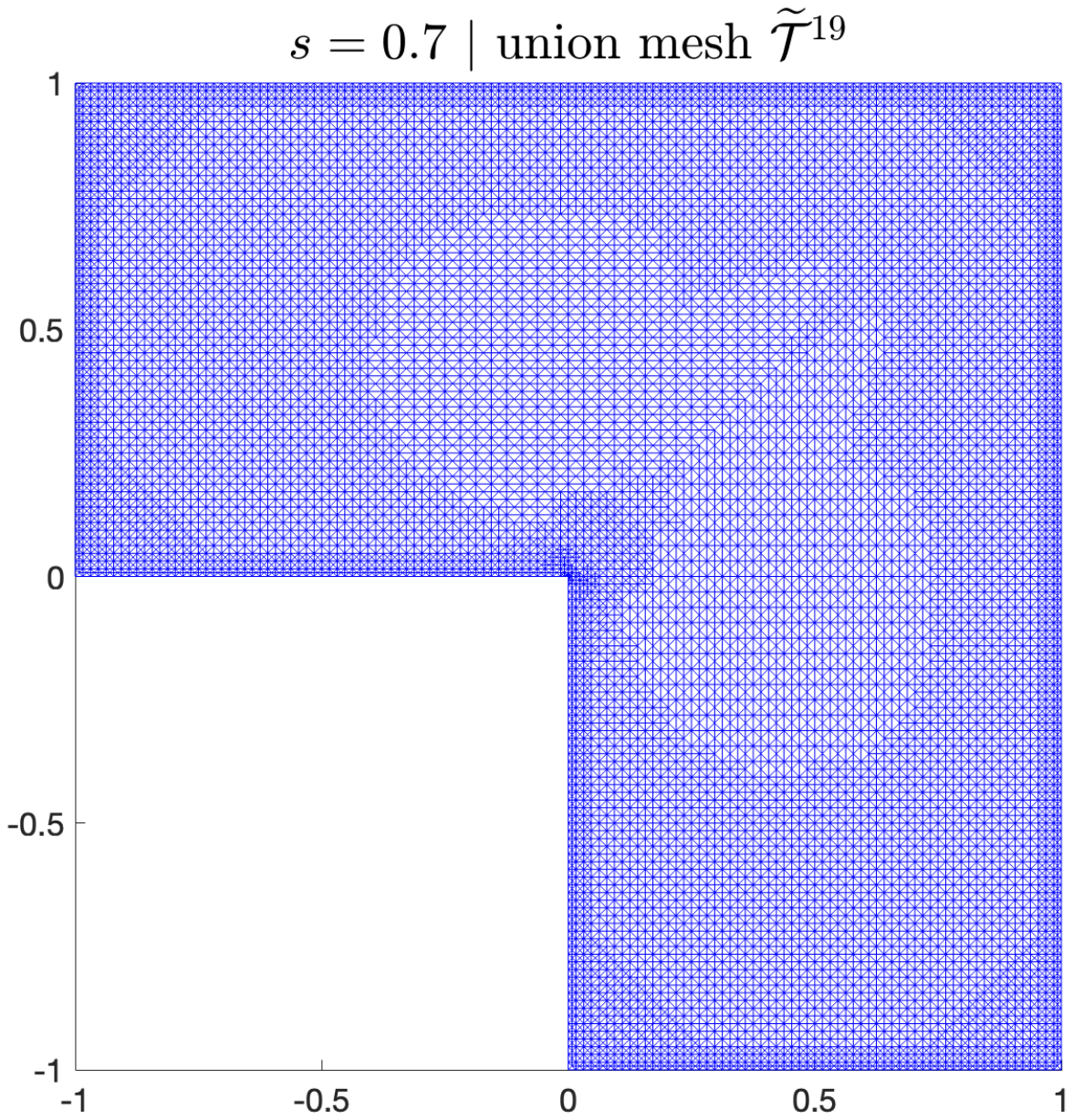}
        \
        \includegraphics*[height=4.5cm, trim = 50 165 19 100, clip]{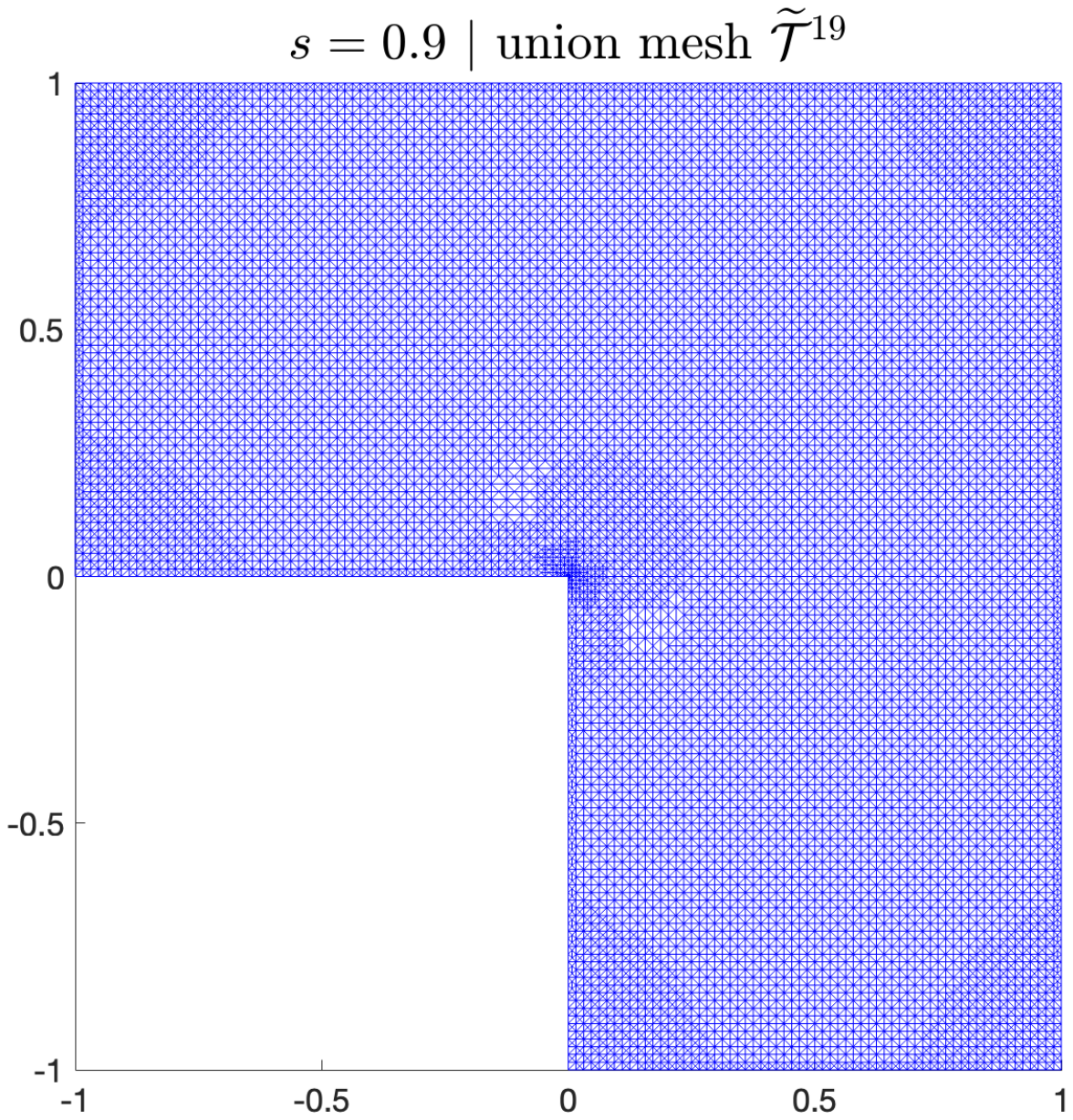}
      \end{center}
    \caption{Test case III (multimesh setting):
    the initial coarse mesh $\T^0$;
    the union meshes $\widetilde\T^{19}$ generated by Algorithm~\ref{alg:algorithm_outline}
    for $s \in \{0.1,\, 0.3,\, 0.5,\, 0.7,\, 0.9\}$.
    }
    \label{fig:test:case:III:meshes}
\end{figure}

In this test case, we use the initial mesh $\T^0$ as a uniform mesh of 384 right-angled triangles
(see Figure~\ref{fig:test:case:III:meshes}).
For each $s \in \{0.1,\, 0.3,\, 0.5,\, 0.7,\, 0.9\}$, we set the same tolerance as
in test case~II (see Table~\ref{tab:test:case:II:rates}) and run our adaptive \emph{multimesh} refinement algorithm.
The decay rates for global a posteriori error estimates $\eta^m$ and $\widetilde\eta^m$
(see~\eqref{eq:triangular_inequality_fractional_estimator} and~\eqref{eq:union_mesh_fractional_estimate_def}, respectively)
as well as the theoretical convergence rates predicted in~\cite[Theorem 4.3]{2015_bonito}
for uniform mesh refinement
are given in Table~\ref{tab:test:case:III:rates}. 

\begin{table}[!th]
  \caption{Test case III (multimesh setting):
  the number $m_*$ of iterations and the decay rates for global error estimates
  $\eta^m$,\, $\widetilde\eta^m$ as well as
  the theoret\-ical convergence rates for uniform mesh
  refinement 
  (see~\cite[Theorem 4.3]{2015_bonito}).
  }
  \label{tab:test:case:III:rates}
  \begin{center}
    \begin{tabular}{lccccc}
\toprule
$s$ & 0.1 & 0.3 & 0.5 & 0.7 & 0.9 \\
\midrule
$m_*$ & 39 & 34 & 33 & 35 & 40\\
\midrule
decay rate for $\eta^m$ & 0.64 & 0.84 & 0.93 & 0.96 & 0.98 \\
decay rate for $\widetilde \eta^m$ & 0.67 & 0.91 & 0.97 & 0.99 & 1.00 \\
decay rate for uniform refinement & 0.35 & 0.55 & 0.67 & 0.67 & 0.67 \\
\bottomrule
\end{tabular}

  \end{center}
\end{table}

As in all previous test cases, we observe that the global error estimates $\eta^m$
based on the triangle inequality decay with a slower rate than the estimates $\widetilde\eta^m$ computed using union meshes.
This conclusion once again confirms superior effectivity and robustness
of the global error estimation strategy based on the union mesh.
In all computations the decay rates for global error estimates $\eta^m$,\, $\widetilde\eta^m$
generated by the adaptive algorithm exceed the rates
predicted by~\cite[Theorem 4.3]{2015_bonito} for uniformly refined meshes.
As expected, for $s = 0.9$, the adaptive algorithm has recovered the optimal convergence rate,
$\mathcal{O}\big(\mathrm{dof}^{-1}\big)$.
We emphasize 
that this optimal rate cannot be achieved 
if uniform mesh refinement is used.
This is due to the presence of dominant geometric singularity in the solution
to problem~\eqref{eq:fractional_laplacian_strong} when $s$ is close to 1.
We note that for $s = 0.7$, the rate is also very close to~optimal.


Figure~\ref{fig:test:case:III:meshes} depicts the initial coarse mesh $\T^0$
as well as the union meshes generated at the 19th iteration of Algorithm~\ref{alg:algorithm_outline}
for $s \in \{0.1,\, 0.3,\, 0.5,\, 0.7,\, 0.9\}$.
This figure shows how adaptive mesh refinement reflects the interplay of two types of singularities
in the solution to problem~\eqref{eq:fractional_laplacian_strong} in this test case for different values of $s$,
with strong refinement solely in the boundary layer for $s=0.1$
gradually transitioning towards
strong refinement at the reentrant corner and barely any refinement in the boundary layer for $s=0.9$.

\section{Conclusions} \label{sec:conclusions}

The numerical solution of FPDEs presents many challenges due to the nonlocal nature of
the underlying pseudo-differential operators.
Adaptive solution strategies have the potential to
substantially reduce computational times and achieve
optimal convergence and complexity in the presence of sharp interfaces,
boundary layers and geometric singularities in solutions to FPDEs.

Focusing on the spectral fractional Laplacian,
an important contribution of this paper is in developing a novel
\emph{multimesh} rational approximation scheme for
discretizing fractional powers of elliptic operators.
The scheme comes with an effective a posteriori error estimation strategy
driving an adaptive multimesh refinement algorithm
that is implemented in two dimensions in the open-source MATLAB package T-IFISS.
Extensive numerical experimentation has demonstrated
the effectivity and robustness of the global error estimation strategy based on the union mesh,
superior convergence rates for approximations generated by the adaptive algorithm
compared to approximations obtained via uniform mesh refinement,
and, crucially, significant reductions in computational complexity
and actual runtimes of the adaptive multimesh refinement algorithm
compared to its single-mesh counterpart.

In this work, adaptive approximations are generated on
\emph{locally quasi-uniform} finite element meshes.
Therefore, the observed convergence rates for adaptive multimesh approximations
are suboptimal in test cases with smaller fractional powers of the Laplacian.
Optimal convergence rates can be recovered by employing \emph{anisotropic} mesh refinement
that captures singular behavior of the solution within boundary/interior layers
more effectively than refinements that preserve shape-regularity of meshes.
An extension of the multimesh approach proposed in this paper to anisotropic meshes,
including a posteriori error estimation on anisotropic elements and, critically,
\rev{an adaptive refinement strategy for anisotropic meshes} will be the subject of future work.

\rev{
Another important extension concerns the design and implementation of a fully-adaptive strategy
that incorporates adaptive selection of the fineness parameter $\kappa$
for the rational approximation of the function $\lambda \mapsto \lambda^{-s}$ (see~\S\ref{sec:rational:approx})
in addition to adaptive finite element approximations
for parametric reaction-diffusion problems (cf.~\cite[Section 7.2]{2023_bulle}
in the single-mesh setting).
In the context of the BP rational approximation scheme~\eqref{eq:bp_rational_approximation},
this involves adding new quadrature points at some iterations of the algorithm
for the numerical approximation of the integral in~\eqref{eq:balakrishnan}.
A critical question that needs to be addressed here is
what finite element meshes one should assign to new parametric problems
arising from newly introduced quadrature points in the rational scheme,
in order to balance the decay of the overall discretization error 
against the growth in the total number of degrees of freedom.
These issues and the associated implementation aspects require further investigation.
}

\section*{Acknowledgments}
This work was initiated when the second author was visiting the School of Mathematics at
the University of Birmingham.
Partial support of that visit by the London Mathematical Society via Scheme~4 (Research in Pairs)
is gratefully acknowledged.
The first author would like to thank the
Institute for Computational and Experimental Research in Mathematics (ICERM) at Brown University
for support and hospitality during the Semester Program
\emph{``Numerical PDEs: Analysis, Algorithms, and Data Challenges''} (January--May 2024),
where part of the work on this paper was undertaken.
The authors also thank Natalia Kopteva (University of Limerick) for stimulating discussions
on the topic of this work and for her helpful suggestions that led to improvements
in the presentation of numerical results.

\bibliographystyle{siam}
\bibliography{biblio.bib}

\begin{thebibliography}{10}

\bibitem{2017_ainsworth}
{\sc M.~Ainsworth and C.~Glusa}, {\em Aspects of an adaptive finite element
  method for the fractional {L}aplacian: a priori and a posteriori error
  estimates, efficient implementation and multigrid solver}, Comput. Methods
  Appl. Mech. Engrg., 327 (2017), pp.~4--35.

\bibitem{2000_ainsworth}
{\sc M.~Ainsworth and J.~T. Oden}, {\em A Posteriori Error Estimation in Finite
  Element Analysis}, Pure Appl. Math. (N. Y.), Wiley, New York, 2000.

\bibitem{2019_antil}
{\sc H.~Antil, Y.~Chen, and A.~Narayan}, {\em Reduced basis methods for
  fractional {L}aplace equations via extension}, SIAM J. Sci. Comput., 41
  (2019), pp.~A3552--A3575.

\bibitem{1960_balakrishnan}
{\sc A.~V. Balakrishnan}, {\em Fractional powers of closed operators and the
  semigroups generated by them}, Pacific J. Math., 10 (1960), pp.~419--437.

\bibitem{BanjaiMNOSS_19_TFS}
{\sc L.~Banjai, J.~M. Melenk, R.~H. Nochetto, E.~Ot\'arola, A.~J. Salgado, and
  C.~Schwab}, {\em Tensor {FEM} for spectral fractional diffusion}, Found.
  Comput. Math., 19 (2019), pp.~901--962.

\bibitem{2023_banjai}
{\sc L.~Banjai, J.~M. Melenk, and C.~Schwab}, {\em Exponential convergence of
  {$hp$} {FEM} for spectral fractional diffusion in polygons}, Numer. Math.,
  153 (2023), pp.~1--47.

\bibitem{1985_bank}
{\sc R.~E. Bank and A.~Weiser}, {\em Some a posteriori error estimators for
  elliptic partial differential equations}, Math. Comp., 44 (1985),
  pp.~283--301.

\bibitem{2021_tifiss_paper}
{\sc A.~Bespalov, L.~Rocchi, and D.~Silvester}, {\em T-{IFISS}: a toolbox for
  adaptive {FEM} computation}, Comput. Math. Appl., 81 (2021), pp.~373--390.

\bibitem{2018_bonito}
{\sc A.~Bonito, J.~P. Borthagaray, R.~H. Nochetto, E.~Ot\'arola, and A.~J.
  Salgado}, {\em Numerical methods for fractional diffusion}, Comput. Vis.
  Sci., 19 (2018), pp.~19--46.

\bibitem{2017_bonitoa}
{\sc A.~Bonito, W.~Lei, and J.~E. Pasciak}, {\em The approximation of parabolic
  equations involving fractional powers of elliptic operators}, J. Comput.
  Appl. Math., 315 (2017), pp.~32--48.

\bibitem{2017_bonito}
\leavevmode\vrule height 2pt depth -1.6pt width 23pt, {\em Numerical
  approximation of space-time fractional parabolic equations}, Comput. Methods
  Appl. Math., 17 (2017), pp.~679--705.

\bibitem{2021_bonito}
{\sc A.~Bonito and M.~Nazarov}, {\em Numerical simulations of surface
  quasi-geostrophic flows on periodic domains}, SIAM J. Sci. Comput., 43
  (2021), pp.~B405--B430.

\bibitem{2015_bonito}
{\sc A.~Bonito and J.~E. Pasciak}, {\em Numerical approximation of fractional
  powers of elliptic operators}, Math. Comp., 84 (2015), pp.~2083--2110.

\bibitem{2016_bonito}
\leavevmode\vrule height 2pt depth -1.6pt width 23pt, {\em Numerical
  approximation of fractional powers of regularly accretive operators}, IMA J.
  Numer. Anal., 37 (2017), pp.~1245--1273.

\bibitem{2023_bulle}
{\sc R.~Bulle, O.~Barrera, S.~P.~A. Bordas, F.~Chouly, and J.~S. Hale}, {\em An
  a posteriori error estimator for the spectral fractional power of the
  {L}aplacian}, Comput. Methods Appl. Mech. Engrg., 407 (2023), pp.~Paper No.
  115943, 27.

\bibitem{2017_cances}
{\sc E.~Canc\`es, G.~Dusson, Y.~Maday, B.~Stamm, and M.~Vohral\'ik}, {\em
  Guaranteed and robust a posteriori bounds for {L}aplace eigenvalues and
  eigenvectors: conforming approximations}, SIAM J. Numer. Anal., 55 (2017),
  pp.~2228--2254.

\bibitem{2021_carstensen}
{\sc C.~Carstensen, A.~Ern, and S.~Puttkammer}, {\em Guaranteed lower bounds on
  eigenvalues of elliptic operators with a hybrid high-order method}, Numer.
  Math., 149 (2021), pp.~273--304.

\bibitem{2015_chen}
{\sc L.~Chen, R.~H. Nochetto, E.~Ot\'arola, and A.~J. Salgado}, {\em A {PDE}
  approach to fractional diffusion: a posteriori error analysis}, J. Comput.
  Phys., 293 (2015), pp.~339--358.

\bibitem{2018_cusimano}
{\sc N.~Cusimano, F.~del Teso, L.~Gerardo-Giorda, and G.~Pagnini}, {\em
  Discretizations of the spectral fractional {L}aplacian on general domains
  with {D}irichlet, {N}eumann, and {R}obin boundary conditions}, SIAM J. Numer.
  Anal., 56 (2018), pp.~1243--1272.

\bibitem{DEliaDGGTZ_20_NMN}
{\sc M.~D'Elia, Q.~Du, C.~Glusa, M.~Gunzburger, X.~Tian, and Z.~Zhou}, {\em
  Numerical methods for nonlocal and fractional models}, Acta Numer., 29
  (2020), pp.~1--124.

\bibitem{1996_dorfler}
{\sc W.~D\"orfler}, {\em A convergent adaptive algorithm for {P}oisson's
  equation}, SIAM J. Numer. Anal., 33 (1996), pp.~1106--1124.

\bibitem{2021_faustmann}
{\sc M.~Faustmann, J.~M. Melenk, and D.~Praetorius}, {\em Quasi-optimal
  convergence rate for an adaptive method for the integral fractional
  {L}aplacian}, Math. Comp., 90 (2021), pp.~1557--1587.

\bibitem{2020_harizanov}
{\sc S.~Harizanov, R.~Lazarov, S.~Margenov, and P.~Marinov}, {\em Numerical
  solution of fractional diffusion-reaction problems based on {BURA}}, Comput.
  Math. Appl., 80 (2020), pp.~316--331.

\bibitem{kpp}
{\sc M.~Karkulik, D.~Pavlicek, and D.~Praetorius}, {\em On 2{D} newest vertex
  bisection: optimality of mesh-closure and {$H^1$}-stability of
  {$L_2$}-projection}, Constr. Approx., 38 (2013), pp.~213--234.

\bibitem{KoptevaR_10_SMN}
{\sc N.~Kopteva and E.~O'Riordan}, {\em Shishkin meshes in the numerical
  solution of singularly perturbed differential equations}, Int. J. Numer.
  Anal. Model., 7 (2010), pp.~393--415.

\bibitem{2020_lischke}
{\sc A.~Lischke, G.~Pang, M.~Gulian, and et~al.}, {\em What is the fractional
  {L}aplacian? {A} comparative review with new results}, J. Comput. Phys., 404
  (2020), pp.~109009, 62.

\bibitem{2015_nochettob}
{\sc R.~H. Nochetto, E.~Ot\'arola, and A.~J. Salgado}, {\em A {PDE} approach to
  fractional diffusion in general domains: a priori error analysis}, Found.
  Comput. Math., 15 (2015), pp.~733--791.

\bibitem{2020_pfeiler}
{\sc C.-M. Pfeiler and D.~Praetorius}, {\em D\"orfler marking with minimal
  cardinality is a linear complexity problem}, Math. Comp., 89 (2020),
  pp.~2735--2752.

\bibitem{tifiss_software}
{\sc D.~J. Silvester, A.~Bespalov, Q.~Liao, and L.~Rocchi}, {\em {T}riangular
  {IFISS} ({T-IFISS})}.
\newblock Available online at \url{http://www.manchester.ac.uk/ifiss/tifiss},
  February 2019.

\bibitem{stevenson}
{\sc R.~Stevenson}, {\em The completion of locally refined simplicial
  partitions created by bisection}, Math. Comp., 77 (2008), pp.~227--241.

\bibitem{2020_vabishchevich}
{\sc P.~N. Vabishchevich}, {\em Approximation of a fractional power of an
  elliptic operator}, Numer. Linear Algebra Appl., 27 (2020), pp.~e2287, 14.

\bibitem{1994_verfurth}
{\sc R.~Verf\"urth}, {\em A posteriori error estimation and adaptive
  mesh-refinement techniques}, J. Comput. Appl. Math., 50 (1994), pp.~67--83.

\bibitem{2017_zhao}
{\sc X.~Zhao, X.~Hu, W.~Cai, and G.~E. Karniadakis}, {\em Adaptive finite
  element method for fractional differential equations using hierarchical
  matrices}, Comput. Methods Appl. Mech. Engrg., 325 (2017), pp.~56--76.

\end{thebibliography}
\end{document}